\documentclass[11pt]{article}

\usepackage[a4paper, top=1.2in, bottom=1.3in,
left=1.2in, right=1.2in]{geometry}
\usepackage[parfill]{parskip}
\usepackage{needspace}
\usepackage{microtype}

\usepackage[runin]{abstract}
\usepackage{fancyhdr}
\usepackage[explicit]{titlesec}
\usepackage[section]{placeins}
\usepackage{multirow}

\usepackage[T1]{fontenc}
\usepackage[cal=esstix,
            scr=boondoxupr,
            frak=euler,
            bb=libus,
            scaled=0.95]{mathalpha}
\usepackage[OT1, small,
            euler-hat-accent,
            euler-digits]{eulervm}
\usepackage{Baskervaldx}
\usepackage{pifont}

\usepackage[leqno]{amsmath}
\usepackage{amssymb}
\usepackage{amsthm}
\usepackage{mathtools,thmtools}
\usepackage{stmaryrd}
\usepackage{xfrac,xparse,xspace}
\usepackage{stackrel,scalerel}

\usepackage[shortlabels]{enumitem}
\usepackage{booktabs,diagbox,makecell,multirow}
\usepackage{arydshln}

\usepackage[english]{babel}
\usepackage{csquotes}
\usepackage[backend=biber,
            style=alphabetic,
            hyperref=auto,
            useprefix=false,
            maxbibnames=99,
            maxalphanames=99,
            maxcitenames=99]{biblatex}  
\usepackage{hyperref}
\usepackage{cleveref}
\usepackage{url}
\usepackage{crossreftools}

\usepackage[dvipsnames]{xcolor}
\usepackage{graphicx}
\usepackage{tikz,tikz-cd,tikz-imagelabels}
\usepackage{dynkin-diagrams}
\usepackage{rotating}
\usepackage{float}
\usepackage[font=small,
            labelfont=bf,
            labelsep=period,
            format=plain,
            singlelinecheck=false]{caption}
\usepackage[all,cmtip]{xy}
 
\usepackage{xcolor}

\newcommand\todo[1]{}
\newcommand\tr{\operatorname{tr}}
\newcommand\Sym{\operatorname{Sym}}
\newcommand\Acon{A_{\operatorname{con}}}
\newcommand\Afib{A_{\operatorname{fib}}}

\newcommand\sheaf{\mathcal{O}}

\newcommand\thestack{\mathfrak{X}}
\newcommand\db{D^b}

\newcommand\Rhom{\mathbf{R}\operatorname{Hom}}
\newcommand\Hom{\operatorname{Hom}}
\newcommand\End{\operatorname{End}}
\newcommand\Ext{\operatorname{Ext}}

\newcommand\CC{\mathbb{C}}
\newcommand\Aff{\mathbb{A}}
\newcommand{\perf}{\mathfrak{Perf}}

\theoremstyle{definition}
\newtheorem{thm}{Theorem}
\numberwithin{thm}{section}
\newtheorem{prop}[thm]{Proposition}
\newtheorem{example}[thm]{Example}

\newtheorem{lemma}[thm]{Lemma}

\theoremstyle{remark}
\newtheorem{remark}[thm]{Remark}

\makeatletter
\newcommand\ackname{Acknowledgements}
\if@titlepage
  \newenvironment{acknowledgements}{%
      \titlepage
      \null\vfil
      \@beginparpenalty\@lowpenalty
      \begin{center}%
        \bfseries \ackname
        \@endparpenalty\@M
      \end{center}}%
     {\par\vfil\null\endtitlepage}
\else
  \newenvironment{acknowledgements}{%
      \if@twocolumn
        \section*{\abstractname}%
      \else
        \small
        \begin{center}%
          {\bfseries \ackname\vspace{-.5em}\vspace{\z@}}%
        \end{center}%
        \quotation
      \fi}
      {\if@twocolumn\else\endquotation\fi}
\fi
\makeatother

\begin{document}

\title{\huge Derived autoequivalences of length $2$ flops via GIT}
\author{Aporva Varshney}
\date{}

\maketitle

\begin{abstract}
   We obtain the derived autoequivalences of a flopping rational curve of length $2$ using GIT and the theory of windows applied to the universal length $2$ flop. 
   We show that the stringy K\"{a}hler moduli space (SKMS) associated to the GIT problem, as constructed by Halpern-Leistner--Sam, matches the description of the space obtained for length $2$ threefolds by Hirano--Wemyss as a quotient of a Bridgeland stability manifold.
   Furthermore, we show that its fundamental group acts via contraction algebra and fibre algebra twists, hence recovering the monodromy action described by Donovan--Wemyss.
   In particular, this shows that the two approaches to building the SKMS agree in this setting.
\end{abstract}

\tableofcontents

\section{Introduction}

\subsection{The stringy K\"{a}hler moduli space}

The stringy K\"{a}hler moduli space (SKMS) of a variety is a conjectural space originating in physics, with links to mirror symmetry where it can be interpreted as a moduli space of complex structures of the mirror manifold.
Consequently, the fundamental group of the SKMS should act on the relevant Fukaya category as autoequivalences. 
Hence, passing through mirror symmetry, this space captures information about the derived symmetries of the variety.

Although these ideas are not mathematically well-defined in full generality, there have been various proposals for how to construct the space in specific settings.
For the case of a single flopping rational curve $C$ in a threefold $X$, Hirano--Wemyss described the space as a quotient of a certain Bridgeland stability manifold \cite{hirano-wemyss}.
This work extended the earlier results of Toda \cite{toda2008stability} which applied to the Atiyah flop and Reid's Pagoda flops, and in particular made rigorous the description of the SKMS given for the former in the physics literature \cite{aspinwall}.

The description of the space turns out to depend only on the \textit{length} invariant of the flop, which is the multiplicity of the rational curve in the scheme theoretic fibre of the contraction, introduced by \cite{katz-morrison}.
The length is some number between $1$ and $6$, and Hirano--Wemyss's description of the SKMS for a length $l$ flop is given by a sphere with punctures at the north and south poles, and $N$ equatorial punctures, where $N$ depends on $l$ as follows.

\begin{center}
    \begin{tabular}{ |p{2cm}||p{1cm}|p{1cm}|p{1cm}|p{1cm}|p{1cm}|p{1cm}|  }
 \hline
 $l$ & 1 & 2 & 3 & 4 & 5 & 6\\
 \hline
 $N$   & 1 & 2 & 4 & 6 & 10 & 12  \\
 \hline
\end{tabular}
\end{center}

\subsection{Monodromy action}

Fixing a basepoint in the SKMS of a flopping space $X$, we wish to describe to autoequivalences of $\db(X)$ coming from loops in the SKMS.
The north and south pole punctures correspond to line bundle twists on $X$ and its flop $X^{\prime}$ respectively, so the more interesting symmetries come from the equatorial punctures.
In the length $1$ case, monodromy around the single puncture corresponds to the `flop-flop' autoequivalence \cite{bridgeland-flops}.
This is the classical spherical twist around $\sheaf_C(-1)$ for the Atiyah flop \cite{st}, while for Reid's Pagoda it is given by Toda's fat spherical twist \cite{toda-fat-spherical}.
In either case, the twist can be written explicitly as follows. Let the base of the flop be
\[
\{x^2 + y^2 = z^2 + t^{2n}\}
\]
so that when $n = 1$ we have the Atiyah flop, and the Pagoda flop when $n > 1$.
There is a sheaf $\mathcal{E} := \mathcal{E}_n$ obtained by taking $\mathcal{E}_1 := \sheaf_C(-1)$ and $\mathcal{E}_{i + 1} := Cone(\mathcal{E}_{i}[-1] \to \sheaf_C(-1))$.
This sheaf is the universal deformation of $\sheaf_C(-1)$ over the algebra $\CC[t]/t^n$.
The autoequivalence can then be written as the cone on the map
\[
\mathbf{R}\Hom(\mathcal{E}, -) \otimes^{\textbf{L}}_{\mathbb{C}[t]/t^n} \mathcal{E} \to \operatorname{Id}.
\]

To generalise this, Donovan--Wemyss considered the non-commutative deformations of $\sheaf_C$ \cite{dw-noncomm-def}. 
They found the associated deformation functor to be representable, giving a sheaf $\mathcal{E}$ deforming $\sheaf_C(-1)$ over a finite-dimensional algebra $\Acon$ called the contraction algebra.  
One can then twist around this sheaf by taking the cone on the map
\[
\mathbf{R}\Hom(\mathcal{E}, -) \otimes^{\textbf{L}}_{\Acon} \mathcal{E} \to \operatorname{Id}.
\]

For all lengths, this gives the monodromy action around the first equatorial puncture.
For successive punctures we can play a similar game, but with certain thickenings of $C$ \cite{dwskms}.
These thickenings define $N$ sheaves $S_i$ supported set-theoretically on $C$, each of which yields an autoequivalence.
There is a set of algebras $A_i$ controlling the deformations of $S_i$, which can be computed explicitly using successive tilting bundles $\mathcal{V}_i$ on $X$.
Each $A_i$ is a quotient of the algebra $\End_{X}(\mathcal{V}_i)$, and we have functors
\[
    \db(A_i) \to \db(\End_{X}(\mathcal{V}_i)) \to \db(X)
\]
where the first map is induced by the quotient map on algebras, and the second map is the derived equivalence induced by the tilting bundle.
These functors allow us to understand the autoequivalences as twists around spherical functors \cite{godinho}.
In particular $A_1 \cong \Acon$, and the twist around the spherical functor $\db(\Acon) \to \db(X)$ agrees with the flop-flop autoequivalence.

For the relevant case of length $2$ flops, the autoequivalence coming from the additional puncture can be constructed using the structure sheaf $\sheaf_{2C}$ of the non-reduced curve.
The algebra $\Afib$ which governs the deformations of this sheaf is called the \textit{fibre algebra}.

\subsection{Universal flops}

Given that the number of punctures $N$ depends only on the length $l$, it is natural to consider the \textit{universal} flop of length $l$, which we denote $\mathcal{X}_l$. 
This is a certain space such that a threefold length $l$ flop can be obtained locally via a suitable classifying map.
The universal flop of length $1$ is the Atiyah flop and Reid's Pagoda flop with base $\{x^2 + y^2 = z^2 - t^{2n}\}$ can be obtained using the base change $s = t^n$ from the conifold $\{x^2 + y^2 = z^2 - s^2\}$.
For length $2$, the universal flop was first studied in \cite{curto-morrison} and is a $6$-dimensional space resolving the singularities of the base
\[
\{x^2 + uy^2 + 2vyz + wz^2 + (uw - v^2)t^2 = 0\}.
\]
For higher lengths, it becomes infeasible to describe the spaces concretely, so instead they were studied algebraically in \cite{Karmazyn2017}, where they were shown to be quiver moduli spaces.

Focusing now on the length $2$ case, one can construct the contraction and fibre algebras for the universal flop using a natural tilting bundle given by \cite{vdb} as in the threefold case.
The work of \cite{dw-noncomm-enhancement} tells us that these algebras similarly induce autoequivalences for the universal flop\footnote{In \textit{loc. cit.} only the contraction algebra is shown to induce an autoequivalence. There are some obstructions to showing that the fibre algebra also induces an autoequivalence in general, but in the specific case of the universal length $2$ flop one can use standard geometric arguments; see \Cref{subsec:fibre}.}
and moreover, as a consequence of \cite{Karmazyn2017}, one can recover the corresponding algebras for the threefolds using the classifying map and the algebras for the universal space.
Hence, one may consider the derived symmetries of length $2$ threefolds as being induced by those of the universal flop.

\subsection{SKMS of a GIT Quotient}

A central motivation for the present work is provided by the alternative proposal for the SKMS in the case of a quasi-symmetric linear reductive GIT quotient $X = \thestack^{ss}/G$, by Halpern-Leistner--Sam \cite{hlsam}.
This definition also has links to the physics literature, and a certain space known as the Fayet–-Iliopoulos parameter space, which is closely linked to the SKMS.

In contrast to the methods outlined above, Halpern-Leistner--Sam instead consider the character lattice of the maximal torus of the group $G$, denoted $M$. 
Within the space $M_{\mathbb{R}} := M \otimes \mathbb{R}$, they build a hyperplane arrangement.
The SKMS is constructed by taking the complement of this arrangement from the Weyl invariant subspace $M_{\mathbb{R}}^{\mathcal{W}}$, complexifying, and finally quotienting by the action of the Picard group.
Following ideas from \cite{hhp, segal}, the monodromy action can then be described through `window-shift autoequivalences,' obtained by finding `window' subcategories $\mathscr{W}_j \subset \db(\thestack)$ which are derived equivalent to any GIT quotient, under a suitably generic choice of linearization.
In particular, fixing such linearizations $l, l^{\prime}$ and denoting the respective quotients by $X, X^{\prime}$, one can then compose equivalences 
\[\db(X) \to \mathscr{W}_j \to \db(X^{\prime}) \to \mathscr{W}_{k} \to \db(X)\]
for choices of windows $\mathscr{W}_j, \mathscr{W}_k$ and obtain autoequivalences of $\db(X)$.

\subsection{Results and structure}

We uniformly study length $2$ flops through the universal flopping space $\mathcal{X}_2$ and show that the two approaches to building the SKMS coincide in both the description of the space and the resulting monodromy action.

Firstly, in \Cref{section:git-problem}, we build $\mathcal{X}_2$ as a linear non-abelian GIT quotient.

\begin{prop}[=\ref{prop:git-problem}]
    The universal flop of length $2$ is a GIT quotient of the stack
    \[
    [V \oplus V^* \oplus (\Sym^2V \otimes \det V^{-1})^{\oplus 2} / GL(V)]
    \]
    where $V$ is a two dimensional vector space.
\end{prop}
In physics terminology, this description puts length $2$ flops in the framework of non-abelian gauged linear sigma models.

The SKMS for this GIT problem is computed in \Cref{section:skms}, and we find that it is precisely given by a sphere with north and south pole punctures and $2$ equatorial punctures.

We then need to check that the monodromy action is through the fibre and contraction algebras to show that the two descriptions of the SKMS align.
For the universal flop, these are given by \cite{kawamata}:
\begin{align*}
    \Afib & \cong \CC[u, v, w] \\
    \Acon & \cong \frac{\mathbb{C}[t]\langle\beta, \gamma\rangle}{\begin{matrix}
        [\beta^2, \gamma], [\gamma^2, \beta], t[\beta, \gamma]
    \end{matrix}}
\end{align*}
where $[-, -]$ denotes the commutator.
In \Cref{section:twists}, we compute the monodromy action of the SKMS and detail the relationship between these autoequivalences and the algebras $\Afib$ and $\Acon$.

\begin{prop}[\S \ref{subsec:fibre}]
    The window-shift autoequivalence corresponding to monodromy around an equatorial puncture in the SKMS is the twist around a spherical functor
    \[
    \db(\Afib) \to \db(\mathcal{X}_2).
    \]
\end{prop}

\begin{prop}[\S \ref{sec:contractiontwist}]
    The window-shift autoequivalence corresponding to monodromy around an equatorial puncture in the SKMS is the twist around a spherical functor
    \[
    \db(\widetilde{A}_{\text{con}}) \to \db(\mathcal{X}_2)
    \]
    where $\widetilde{A}_{\text{con}}$ is the upper triangular algebra
    \[
    \begin{bmatrix}
        \Acon/[\beta, \gamma] & \Acon/(t, [\beta, \gamma]) \\
        0 & \Acon/t
    \end{bmatrix}
    \]
    formed by gluing the $\Acon$ modules $\Acon / t$ and $\Acon / [\beta, \gamma]$ along the bimodule $\Acon / (t, [\beta, \gamma])$.
\end{prop}
The algebra $\widetilde{A}_{\text{con}}$ appears naturally in the process of deriving the monodromy action from the GIT problem.
To explain its relation to $\Acon$, we draw attention to the case of the node $B := \CC[x, y]/xy$, which has a non-commutative resolution given by the upper triangular algebra
\[
\begin{bmatrix}
    B/(x)  & B/(x, y) \\ 0 & B/(y)
\end{bmatrix}.
\]
Our algebra $\widetilde{A}_{\text{con}}$ similarly plays the role of a non-commutative resolution for $\Acon$, which can be considered as a gluing of two branches in light of the relation $t[\beta, \gamma] = 0$.
This is summarised in the following.

\begin{prop}[=\ref{prop:smooth-algebra}, \ref{prop:equivalent-twists}]
    The algebra $\widetilde{A}_{\text{con}}$ has finite global dimension and admits a pair of adjoint functors $\Phi^L \dashv \Phi$
\begin{align*}
    & \Phi : \db(\widetilde{A}_{\text{con}}) \to \db(\Acon) & \Phi^L : \perf(\Acon) \to \db(\widetilde{A}_{\text{con}})
\end{align*}
such that $\Phi\Phi^L$ is the identity on $\perf(\Acon)$.
\end{prop}
As a consequence, we find that the spherical twists around our functor $\db(\widetilde{A}_{\text{con}}) \to \db(X)$ and the Donovan--Wemyss functor $\db(\Acon) \to \db(X)$ are naturally isomorphic.

We draw attention to the fact our algebra $\tilde{A}_{\text{con}}$ arises naturally from our GIT problem.
One should be able to construct the category $\db(\Acon)$ by quotienting $\db(\tilde{A}_{\text{con}})$ by the kernel of the spherical functor $\db(\widetilde{A}_{\text{con}}) \to \db(\mathcal{X}_2)$.
This kernel is comprised of sheaves on the stack which become zero upon restriction to the semistable locus.
In this way, it would be possible to find the contraction algebra by considering the spherical functors constructed, and thus, our GIT approach gives an alternate derivation of the contraction algebra.

More generally, we can ask whether the autoequivalences and the description of the SKMS obtained in \cite{hirano-wemyss, dwskms} for a single curve flop of \textit{any} length can be obtained using the procedures outlined above applied to universal flops of higher length.
Additionally, it would be interesting to check whether the sequence of tilting bundles constructed for the threefolds by Donovan--Wemyss can also be obtained using the GIT machinery in the length $\geq 3$ cases, and we hope the results here provide some insight in this direction.

\subsection*{Conventions}
We work over $\CC$. The word `module' refers to right modules unless otherwise stated, and for a $\CC$-algebra $A$ the category $\db(A)$ refers to the bounded derived category of finitely generated right modules while for a scheme $X$ the category $\db(X)$ refers to the bounded derived category of coherent sheaves.
For arrows $f, g$ in a quiver, the notation $fg$ refers to the composition $g \circ f$.
For a connected reductive group $G$ with Borel subgroup $B$ and maximal torus $T$, $M$ denotes the character lattice of $T$ and $N$ the cocharacter lattice.
For a lattice $K$, and $K_{\mathbb{R}} := K \otimes \mathbb{R}$ with $K_{\mathbb{\CC}}$ defined analagously.
The notation $M_{\mathbb{R}}^+, M_{\mathbb{R}}^-$ refers to the dominant and anti-dominant chambers of $M_{\mathbb{R}}$, and similarly for $N_{\mathbb{R}}$.
The Weyl group of $G$ is denoted $\mathcal{W}$ and $(-)^{\mathcal{W}}$ denotes taking Weyl-invariants.
For a map between schemes $f : X \to Y$, the functors $f_*$ and $f^*$ are implicilty derived, as are tensor products.
We denote $\End(A) := \operatorname{Ext}^{\bullet}(A, A)$.
For a functor $\Phi$ its left and right adjoints are $\Phi^L$ and $\Phi^R$ respectively.

\begin{acknowledgements}
I am hugely indebted to my advisor Ed Segal for introducing me to this topic and for continued patient guidance. I am also grateful to Michael Wemyss, Will Donovan and Marina Godinho for taking the time to explain various aspects of their work, to Michela Barbieri and Calum Crossley for various discussions, and to Parth Shimpi for conversations regarding \Cref{remark:hyperplanes}.
\end{acknowledgements}

\section{Preliminaries}

\subsection{Contraction and fibre algebras}

We work in the context of a projective birational morphism $f : X \to Y := \operatorname{Spec} R$ where $X$ has at worst Gorenstein terminal singularities, and where $Y$ is normal. Additionally, we assume that $\mathbf{R}f_{*}\sheaf_X \cong \sheaf_{Y}$, where the fibres have dimension at most $1$.

We review the construction of contraction and fibre algebras, following \cite{dw-noncomm-def, dw-contractions-deformations}; for our purposes, the deformation-theoretic interpretations will be less relevant and so we focus on their relation to derived autoequivalences.
By \cite{vdb}, there exists a vector bundle $V$ on $X$ so that $\sheaf \oplus V$ is a tilting bundle, and in particular there are equivalences:
\begin{align*}
    & \Lambda_0 := \End_X(\sheaf \oplus V) \cong \End_{R}(R \oplus f_{*}V) \\
    & \db(X) \simeq \db(\Lambda_0).
\end{align*}
Similarly we also have equivalences
\begin{align*}
    & \Lambda_{-1} := \End_X(\sheaf \oplus V^*) \\
    & \db(X) \simeq \db(\Lambda_{-1}).
\end{align*}
For a summand $A$ of a vector bundle $A \oplus B$, we denote by $[A]$ the two sided ideal of $\End(A \oplus B)$ of morphisms factoring through the set
\[
\operatorname{add}(A) := \{M \;|\; M \; \text{is a summand of} \; A^{\oplus j}, \; j \; \text{finite}\}.
\]
Then the contraction and fibre algebras are defined by
\begin{align*}
\Acon := \Lambda_0 / [\sheaf] &&    \Afib := \Lambda_{-1} / [V^*].
\end{align*}

Now assume that $\dim X = 3$ and the exceptional locus is irreducible.
Let $\mathcal{E}_{\text{con}}$ be the object given by mapping $\Acon$ under the derived equivalence $\db(\Lambda_0) \to \db(X)$.
Then there is an induced autoequivalence given by the cone on
\[
\Rhom(\mathcal{E}_{\text{con}}, -) \otimes_{\Acon} \mathcal{E}_{\text{con}} \to Id.
\]
The same holds for $\Afib$, using the derived equivalence $\db(\Lambda_{-1}) \to \db(X)$.

In \cite{godinho}, it is shown that the autoequivalences of flopping threefolds can be written as twists around spherical functors\footnote{
    In \cite{godinho} only the case of the contraction algebra is detailed. However, the methodology extends to $\Afib$, as well as the other deformation algebras in higher length cases: the key is that since all the algebras involved are local in the threefold setting and finite over a Notherian ring, they are self-injective. We thank Marina Godinho for explaining this to us.
}.
In particular, for a length $2$ threefold the functors given by the compositions
\begin{align*}
    & \db(\Acon) \to \db(\Lambda_0) \to \db(X) \\
    & \db(\Afib) \to \db(\Lambda_{-1}) \to \db(X)
\end{align*}
are spherical.
Note that when discussing spherical functors, one can run into foundational issues; in our context, we overcome these as we are only ever working with derived categories of algebras and varieties.

\subsection{The universal flop of length 2}

The universal flop of length $2$, which we denote by $\mathcal{X}_2$, can be obtained, as in \cite{curto-morrison}, by performing a Grassmann blowup of the singular hypersurface
\[
\operatorname{Spec} R := \{x^2 + uy^2 + 2vyz +wz^2 +(uw - v^2)t^2 = 0\} \subset \mathbb{A}^7.
\]
We give an overview of the geometry of this space and its resolution; parts of this description were also given in \cite{kawamata}.
The singular locus is determined by the following equations.
\begin{align*}
    & x = 0 && uy + vz = 0 \\
    & vy + wz = 0 && z^2 + ut^2 = 0 \\
    & y^2 + wt^2 = 0 && yz - vt^2 = 0 \\
    & (uw - v^2)t = 0 &&
\end{align*}
This locus is reducible, with one component $Z_1$ given by the locus where $\{x = y = z = t = 0\}$, so that it is isomorphic to $\mathbb{A}^3$ with coordinates $(u, v, w)$. 
The second component $Z_2$ is slightly harder to describe, but we note that in the open $\{t \neq 0 \}$ it is given by the equations
\begin{align*}
    u & = -z^2t^{-2} \\
    w & = -y^2 t^{-2} \\
    v & = yzt^{-2}
\end{align*}
so that in this open the component is isomorphic to $\mathbb{A}^2_{y, z} \times \mathbb{C}^*_{t}$.
Then, $Z_2$ is given by the closure of this locus, and intersects with $Z_1$ in the locus $\{uw = v^2\}$.
One sees that $Z_2$ can be described as the image of $\mathbb{A}_{t, b, c}^3 \to \mathbb{A}^7$ under the map 
\[t = t, u = b^2, w = c^2, v = bc, x = 0, y = bt, z = ct.\]
We note that generically where $t \neq 0$ this map is an isomorphism since for each fixed $t$ it maps $\mathbb{A}_{b ,c}^2 \to \mathbb{A}^2_{y, z}$, but at $t = 0$ it is the quotient map $\mathbb{A}^2 \to \{uw = v^2\} \subset \mathbb{A}^3$.

Resolving the singularities, we obtain our universal flop.
The exceptional locus is also reducible, with components $E_1 \to Z_1$ and $E_2 \to Z_2$.
The component $E_1$ can be described as a family of conics over $Z_1$:
\[
E_1 \cong \{u\phi^2 + v\phi\epsilon + w\epsilon^2 + \delta^2 = 0\} \subset \mathbb{P}^2_{[\phi : \epsilon : \delta]} \times \mathbb{A}^3_{u, v, w}
\]
so that generically over the locus $\{uw = v^2\}$ the conics degenerate and become reducible, and at the origin we have a non-reduced $\mathbb{P}^1$.
Note that the locus where the conics degenerate is the intersection of $E_1$ and $E_2$.
Outside of this intersection, the generic fiber of $E_2$ is a $\mathbb{P}^1$.
Furthermore, $E_2$ can be described through the diagram
\[\begin{tikzcd}
	{\tilde{E}_2 \cong \mathbb{A}_3\times \mathbb{P}^1} && {E_2} \\
	\\
	{\mathbb{A}_3} && {Z_2}
	\arrow[from=1-1, to=1-3]
	\arrow[from=1-1, to=3-1]
	\arrow[from=1-3, to=3-3]
	\arrow[from=3-1, to=3-3]
\end{tikzcd}\]
where the map $\mathbb{A}_3 \to Z_2$ is as above, and the map $\tilde{E}_2 \to E_2$ is generically an isomorphism. 
At a generic point of the locus $\{t = 0\} \subset Z_2$, the diagram becomes the normalization of a degenerate conic
\[\begin{tikzcd}
	{\mathbb{P}^1 \sqcup \mathbb{P}^1} && {\text{degenerate conic}} \\
	\\
	{\text{pt} \sqcup \text{pt}} && {\text{pt} \in \{uw = v^2\}}
	\arrow[from=1-1, to=1-3]
	\arrow[from=1-1, to=3-1]
	\arrow[from=1-3, to=3-3]
	\arrow[from=3-1, to=3-3]
\end{tikzcd}\]
and at the origin, the map $\tilde{E}_2 \to E_2$ is set-theoretically bijective but scheme-theoretically a double cover. 
Details of this construction are briefly covered in \Cref{prop:factor-acon}, but note that we can also construct the map $\tilde{E}_2 \to E_2$ by blowing up the space $\mathcal{X}_2$ along $E_1$.

The universal flop of length $2$ also arises as quiver GIT applied to a certain quiver with relations.
The details of how the quiver is derived can be found in \cite{Karmazyn2017}, but as a sketch, the basic idea of the construction is to start with the quiver which has as its path algebra the preprojective algebra associated to the Dynkin diagram of type $D_4$. 
We then add deformation parameters; this gives us the quiver for the versal deformation of the minimal resolution of the $D_4$ singularity. 
Since we require the versal deformation of a partial resolution, we have to finally collapse the vertices corresponding to unresolved curves.
This results in the following algebra.

\begin{align}
\begin{aligned}
\begin{tikzpicture}\label{diagram:quiver} [bend angle=15, looseness=1]
\node (C0) at (0,0)  {$0$};
\node (C1) at (2,0)  {$1$};
\draw [->,bend left] (C0) to node[above]  {\scriptsize{$\alpha$}} (C1);
\draw [->,bend left] (C1) to node[below]  {\scriptsize{$\alpha^*$}} (C0);
\draw [->, looseness=24, in=52, out=128,loop] (C1) to node[above] {$\scriptstyle{\beta}$} (C1);
\draw [->, looseness=24, in=-38, out=38,loop] (C1) to node[right] {$\scriptstyle{\gamma}$} (C1);
\draw [->, looseness=24, in=-128, out=-52,loop] (C1) to node[below] {$\scriptstyle{\delta}$} (C1);
\end{tikzpicture}
\end{aligned}
& \quad \quad
\begin{aligned}
&\alpha \alpha^*=te_0, & & \beta^2=T^\beta_0 e_1, \\ 
&\gamma^2=T^\gamma_0 e_1, & & \delta^2=T^\delta_0 e_1, \\ 
&\alpha^*\alpha+\beta+\gamma + \delta = \frac{t}{2} \, e_1 \\ 
\end{aligned}
\end{align}
Here $e_0$ and $e_1$ denote the identities at the respective vertices.
By following King's construction \cite{king}, we can recover the universal length $2$ flop as a moduli space of semistable representations with dimension vector $(1, 2)$.

From this quiver, we can also compute the fibre and contraction algebras, by quotienting by the maps factoring through the second vertex and the first vertex respectively.
The computation of the fibre algebra is immediate: all maps in the quiver are zero, as well as the endomorphism $t$ at the vertex since $\alpha = \alpha^* = 0$, so we are left with just the remaining deformation parameters.
\begin{align*}
    \Afib \cong \mathbb{C}[T_0^\beta, T_0^\gamma, T_0^\delta].
\end{align*}

Simplifying the presentation of the contraction algebra is slightly more difficult, but one can show that \cite[Theorem 4.5]{kawamata}:
\begin{align*}
    \Acon \cong \frac{\mathbb{C}[t]\langle\beta, \gamma\rangle}{\begin{matrix}
        [\beta^2, \gamma], [\gamma^2, \beta]\\ t[\beta, \gamma]
    \end{matrix}}.
\end{align*}
These relations tell us that $\beta^2, \gamma^2$ are central elements and that $\beta, \gamma$ commute after multiplication by $t$.

We consider the following observations, which were also stated in \cite[Prop. 4.9]{kawamata}.
Firstly, using the change of basis given in \cite[Eq. 5.1]{Karmazyn2017}, when $x = y = z = t = 0$ we have that the deformation parameters are linearly related to the coordinates $u, w, v$ as follows.
\begin{align*}
    & u = -T_0^\beta  &&
w = -T_0^\gamma &&
v = -(T_0^\beta + T_0^\gamma - T_0^{\delta})/2
\end{align*}
Therefore, the coordinate ring of $Z_1$ is isomorphic to $\Afib$ as an $R$-module.

Secondly, using \cite{Karmazyn2017} again, the structure of $\Acon$ as an $R$-module is given through the following equations.
\begin{align*}
    & u = \beta^2  && w = \gamma^2  && v = \frac{1}{2}(\beta \gamma + \gamma \beta)\\
    & z = -t\beta && y = -t\gamma  && t = t \\
    & x = 0
\end{align*}
We see that $\Acon$ is supported on $Z_1 \cup Z_2$, as expected by the results of \cite{dw-noncomm-enhancement}.
The abelianisation of $\Acon$ is given by $\mathbb{C}[t, b, c]$, and as an $R$-module this is then supported only on $Z_2$: note that the map $R \to \Acon$ induces a map $R \to \mathbb{C}[t, b, c]$ which precisely agrees with the description of $Z_2$ given above.

We end with an example of how this picture relates to threefold flops.

\begin{example}
    Consider the Laufer flop \cite{laufer}, which has base with equation
    \[
        \{ x^2 + y^3 - tz^2 - yt^3 = 0 \}.
    \]
    From the universal flop, this can be found by the slicing 
    \[
        w = -t, u = y, v = 0.
    \]
    Using the $R$-module structure on $\Acon$ given above, we find that
    \[
    \Acon \otimes_{R} R/(w+t, u-y, v) \cong \mathbb{C}\langle \beta, \gamma \rangle / (\beta^2 - \gamma^3, \beta\gamma + \gamma\beta)
    \]
    which is indeed the contraction algebra for the Laufer flop \cite{dw-noncomm-def}.
\end{example}

\section{Quasi-symmetric linear GIT problem}\label{section:git-problem}

Our first step is to show that we can obtain the universal flop through a linear, quasi-symmetric quotient by $GL_2$.

Beginning with quiver GIT applied to diagram \eqref{diagram:quiver}, there are two choices of stability which yield a good moduli space, corresponding to the two sides of the flop. 
The first stability condition $\theta_1$ dictates that the image of $\alpha$ in $\mathbb{C}^2$ is non-zero, and that at least one of $\beta, \gamma$ and $\delta$ do not preserve this line. 
The other stability condition $\theta_2$ dictates that $\alpha^*$ is non-zero, and that at least one of $\beta, \gamma$ and $\delta$ do not preserve the line given by $\ker \alpha^*$.

\begin{lemma}\label{lemma:tracefree}
    In any $\theta_1$ semistable quiver representation, the maps $\beta, \gamma, \delta$ are trace-free.
\end{lemma}
\begin{proof}
    Under $\theta_1$, we may assume that $\alpha$ is non-zero. 
    From the relations on the quiver \eqref{diagram:quiver}, we see that each of the maps $\beta, \gamma$ and $\delta$ is either trace-free or some multiple of the identity.
    
    The relation
    \[
    \alpha^*\alpha + \beta + \gamma + \delta = \frac{t}{2}I_{2}
    \]
    implies, by taking traces, that $\tr{\beta} + \tr{\gamma} + \tr{\delta} = 0$, hence if any two of the maps are trace-free then so is the third.
    Assume without loss of generality that $\beta$ and $\gamma$ are not trace-free: since they are then multiples of the identity, they both preserve $\alpha$, and then $\delta$ is forced to be trace-free.
    Now we deduce that $\beta + \gamma = 0$, and find that
    \begin{align*}
        & \alpha^*\alpha + \delta = \frac{t}{2}I_2 \implies \alpha\alpha^*\alpha + \alpha\delta = \frac{t}{2}\alpha \implies \alpha\delta = -\frac{t}{2}\alpha.
    \end{align*}
    This contradicts stability, so $\beta, \gamma, \delta$ are indeed trace-free in any semistable representation.
\end{proof}

Fixing our stability condition to be $\theta_1$, using \Cref{lemma:tracefree}, we may now assume $\beta, \gamma, \delta$ are trace-free.

\begin{prop}\label{prop:git-problem}
    The universal flop of length $2$ is a GIT quotient of the stack
    \[
    \thestack := [V \oplus V^* \oplus (\Sym^2V(-1))^{\oplus 2} / GL(V)]
    \]
    where $V$ is a two dimensional vector space and we denote by $\Sym^2V(-1) := \Sym^2V \otimes \det V^{-1}$ the representation of trace free endomorphisms of $V$. 
\end{prop}
\begin{proof}
Consider the GIT problem built via quiver GIT.
This is an affine space with coordinates coming from the maps in the quiver:
\begin{align*}
    \begin{aligned}
    & \alpha = \begin{bmatrix}
           \alpha_{0} & \alpha_{1}
         \end{bmatrix} & 
    & \alpha^* = \begin{bmatrix}
           \alpha^{*}_{0} \\ \alpha^{*}_{1}
         \end{bmatrix} & 
    & \beta = \begin{bmatrix}
            b_{00} & b_{01} \\ b_{10} & -b_{00}
         \end{bmatrix} & 
    & \gamma = \begin{bmatrix}
            c_{00} & c_{01} \\ c_{10} & -c_{00}
         \end{bmatrix} &
    & \delta = \begin{bmatrix}
            d_{00} & d_{01} \\ d_{10} & -d_{00}
         \end{bmatrix}
\end{aligned}
\end{align*}
and additionally the deformation parameters $t, T_0^{\beta}, T_0^{\gamma}, T_0^{\delta}$.
We see that the deformation parameters can be removed since there are relations
\begin{align*}
    t & = \alpha_0\alpha_0^* + \alpha_1\alpha_1^* \\
    T_0^{\beta} & = b_{00}^{2} + b_{01}b_{10} \\
    T_0^{\gamma} & = c_{00}^{2} + c_{01}d_{10} \\
    T_0^{\delta} & = d_{00}^{2} + d_{01}d_{10}
\end{align*}
and additionally we can remove $\delta$ using the relation $\delta = \tfrac{t}{2} - (\beta + \gamma + \alpha^*\alpha)$.
So finally we are left with an affine space with coordinates coming from the maps $\alpha, \alpha^*, \beta, \gamma$.
Denoting by $V \cong \mathbb{C}^2$ the two dimensional vector space at vertex $1$, we see that there is a natural action of $\mathbb{C}^* \times GL(V)$.
Quotienting out by the diagonal $\mathbb{C}^*$, we are left with only a $GL(V)$ action. 
As representations of $GL(V)$, the affine space is precisely of the form claimed in the statement of the proposition.
\end{proof}

Recall that a finite dimensional representation $W$ of a connected reductive group $G$ is quasi-symmetric \cite[\S 1.6]{svdb-nccr} if for every line $l \subset M_{\mathbb{R}}$ through the origin, the subset of weights $\beta \in \Psi$ contained in $l$ sum to zero, where $\Psi$ is the set of weights of the representation.
In particular, any self-dual representation is quasi-symmetric. 
Furthermore, for our problem the generic stabilizer under the action of the maximal torus is finite, hence we satisfy the hypotheses to apply the methods of \cite{hlsam, svdb-nccr, schobers} to our GIT quotient.

\section{Windows and the stringy K\"{a}hler moduli space}\label{section:skms}

We recall the construction of the SKMS following \cite{hlsam}.
The first step is to define a Weyl invariant polytope $\overline{\nabla} \subset M_{\mathbb{R}}$:
\[
\overline{\nabla} = \left\{ \chi \in M_{\mathbb{R}} \; \Big\vert \; |\langle \lambda, \chi \rangle| \leq \frac{\eta_{\lambda}}{2} \;\; \text{for all} \;\; \lambda : \mathbb{G}_m \to T \right\}.
\]
Here, the number $\eta_{\lambda}$ associated to a one-parameter subgroup $\lambda$ can be computed combinatorially \cite{hl, hl-kemboi}:
\[
\eta_{\lambda} = \sum_{\alpha \in \Phi} \min(0, \langle \lambda, \alpha\rangle) - \sum_{\beta \in \Psi} \min(0, \langle \lambda, \beta\rangle)
\]
where $\Phi$ is the set of roots of $G$ and $\Psi$ is the set of weights of the representation under the maximal torus $T \subset G$.

\begin{figure}[!h]
    \centering
    \includegraphics[width=0.5\linewidth]{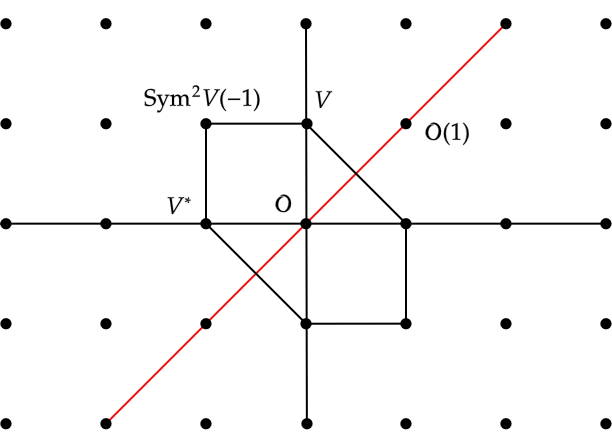}
    \caption{The polytope $\overline{\nabla}$ in the character lattice of the maximal torus, with key representations placed on their highest weights. Red line denotes the Weyl invariant line.}
    \label{fig:zonotope-origin}
\end{figure}

In our setting, the roots of $G = GL_2$ are given by the vectors $\{ (1, -1), (-1, 1)\}$ and the representation has weights 
\begin{align*}
    & V && (1, 0), (0, 1) \\
    & V^* && (-1, 0), (0, -1) \\
    & \Sym^2V(-1) && (1, -1), (0, 0), (-1, 1)
\end{align*}
with the appropriate multiplicities.
We then compute that for a one parameter subgroup 
\[
\lambda(t) = \begin{bmatrix}
    t^m & 0 \\ 0 & t^n
\end{bmatrix}
\]
that the resulting inequality is given by
\[
| mx + ny| \leq \tfrac{1}{2}(|m| +|n| + |m - n|). 
\]
One sees that only $3$ inequalities are needed:
\begin{align*}
    & |x| \leq 1 \\
    & |y| \leq 1 \\
    & |x + y| \leq 1
\end{align*}
which gives us the polytope $\overline{\nabla}$ in \Cref{fig:zonotope-origin}.

To define the SKMS, Halpern-Leistner--Sam construct a hyperplane arrangement in $M_{\mathbb{R}}$ as follows.
Let $(H_i)_{i}$ denote the set of hyperplanes in $M_{\mathbb{R}}$ such that $H_i \cap \overline{\nabla}$ is a facet of $\overline{\nabla}$.
Define a hyperplane arrangement by
\[
\mathcal{H} := \bigcup_{i} \; \{ m - H_{i} \; \vert \; m \in M\}.
\]

\begin{figure}[!h]
    \centering
    \includegraphics[width=0.3\linewidth]{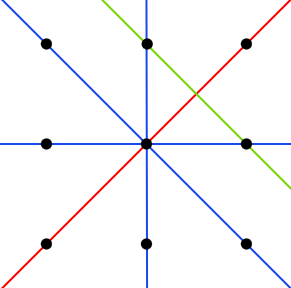}
    \caption{The resulting hyperplane arrangement in blue and green with translates omitted, and Weyl invariant line shown in red.}
    \label{fig:hyerplane-arrangement}
\end{figure}

In our context, the hyperplane arrangement coming from translates of the facets of $\overline{\nabla}$ is partially visualised, omitting translations, in \Cref{fig:hyerplane-arrangement}. 
Using this, the description of the SKMS is given by:
\[
\left( M^{\mathcal{W}}_{\CC} \setminus \bigcup_{H \in \mathcal{H}} (H \cap M^{\mathcal{W}}_{\mathbb{R}}) \otimes \mathbb{C}  \right) \bigg/ M^{\mathcal{W}}.
\]

\begin{prop}
    The SKMS associated to the quasi-symmetric GIT problem 
    \[
    [V \oplus V^* \oplus (\Sym^2V(-1))^{\oplus 2} / GL(V)]
    \]
    is given by a sphere with north and south pole punctures, and two equatorial punctures.
\end{prop}
\begin{proof} 
    Take the obvious identification of the Weyl invariant line $M^{\mathcal{W}}_{\mathbb{R}}$ with $\mathbb{R}$ given by $(t, t) \mapsto t$.
    To obtain the SKMS, we take the complement of the hyperplane arrangement which punctures the Weyl invariant line at integer points and at half integer points. 
    After complexifying and quotienting by the action of $M^{\mathcal{W}}$, which acts via translating by an integer, we obtain a sphere with north and south pole punctures and two equatorial punctures.
\end{proof}

The intersection of the hyperplane arrangement $\mathcal{H}$ with $M^{\mathcal{W}}_{\mathbb{R}}$ defines a poset of faces $\mathcal{C}$ of $M^{\mathcal{W}}_{\mathbb{R}}$.
In our case, we identify $M^{\mathcal{W}}_{\mathbb{R}}$ with $\mathbb{R}$ using $(t, t) \to t$.
Then the hyperplane arrangement in $M_{\mathbb{R}}$ gives us zero-dimensional faces of $M^{\mathcal{W}}_{\mathbb{R}}$ consisting of the points $D_{j} := \tfrac{1}{2}(j + 1)$ and one-dimensional faces consisting of the intervals $C_{j} := (\tfrac{1}{2}j, \tfrac{1}{2}(j + 1))$.
Note that $\overline{C}_{j} \cap \overline{C}_{j + 1} = D_{j}$.

For any face $C \in \mathcal{C}$ and $\delta \in C$, the category split-generated by the characters in $\delta + \overline{\nabla}$ considered as bundles on $\thestack$ does not depend on the choice of $\delta \in C$ \cite[Lemma 5.3]{schobers}, so we denote it by $\mathscr{E}_{C}$, and by $\overline{\nabla}_{C}$ the corresponding polytope.

\begin{thm}[{\cite[Theorem 3.2]{hlsam}}]
    Fix characters $\delta, l \in M^{\mathcal{W}}_{\mathbb{R}}$, with $\delta \in C$ for some face $C$ such that $\partial(\delta + \overline{\nabla}) \cap M = \emptyset$ and $\thestack^{ss}(l) = \thestack^s(l)$.
    Then there is an equivalence given by the restriction functor
    \[
    \db(\thestack^{ss}(l)/G) \simeq \mathscr{E}_C.
    \] 
\end{thm}

For faces $C$ such that the above equivalence holds, we call the subcategories $\mathscr{E}_C$ \textit{windows}.

In our setting, for a one-dimensional face $C_{j}$, the polytope $\overline{\nabla}_{C_{j}}$ contains three lattice points.
The corresponding subcategory $\mathscr{E}_{C_{j}}$ is then a window subcategory, and we denote it by $\mathscr{W}_j := \mathscr{E}_{C_{j}}$ to slightly ease notation.
For a given $j$, the window is given by:
\begin{align*}
    \mathscr{W}_{j} := \begin{cases}
        \langle \sheaf(i), V(i)\rangle  & j = 2i \\
    \langle \sheaf(i), V(i - 1) \rangle & j = 2i - 1
    \end{cases}.
\end{align*}
Note that it follows from \cite[Corollary 1.2.1]{karmazyn2017quiver} that the restriction of $\sheaf \oplus V$ to the GIT quotient recovers the tilting bundle given by Van den Bergh \cite{vdb}.

For use later, we also note that the larger subcategories $\mathscr{E}_{D_j}$ are given as follows.
\[
\mathscr{E}_{D_j} =
    \begin{cases}
        \langle \sheaf(i), V(i), \sheaf(i + 1) \rangle  & j = 2i \\
        \langle \sheaf(i), V(i), V(i - 1), \Sym^2V(i - 1) \rangle & j = 2i - 1
    \end{cases}
\]

\begin{remark}\label{remark:hyperplanes}
    It is natural to ask whether the hyperplane arrangement constructed here has any relation to the hyperplane arrangement constructed in \cite{iyama-wemyss} and used in \cite{hirano-wemyss} to build the SKMS. 
    In fact, the relevant $D_4$ case is explained in \cite[Example 2.9]{iyama-wemyss}: here an infinite set of hyperplanes through the origin in $\mathbb{R}^2$ is used, and an infinite hyperplane arrangement in $\mathbb{R}$ is produced by taking a height $1$ slice.
    Although the method we use produces the same resulting hyperplane arrangement in $\mathbb{R}$, the preceding steps appear to be vastly different.
    For one, our hyperplanes are affine.
    Furthermore, the ambient $\mathbb{R}^2$ comes from different sources; in \cite{iyama-wemyss}, the ambient space is $2$ dimensional since there is a single curve being flopped, while our ambient space is $2$ dimensional due to the \textit{length} of the flop.
    We also point out the dependence of our hyperplane arrangement on the specific GIT problem; in theory, there may be other presentations which yield the same final answer but use a different initial arrangement.
    This leads us to believe that any possible connection between the two approaches must be much deeper.
\end{remark}

\section{Autoequivalences}\label{section:twists}

\subsection{Constructing spherical functors}\label{sec:algorithm}

Here we recall the procedure used to find the monodromy action of the SKMS via autoequivalences. 
We limit the scope of the presentation slightly to reflect our setting.

\subsubsection{Window shift autoequivalences} 
The stability conditions $\theta_1, \theta_2$ on the original quiver induce stability conditions on the new GIT problem. 
Denote by $X_+$ the quotient under $\theta_2$, and $X_-$ the quotient under $\theta_1$.
One may verify that $X_{\pm}$ is the GIT quotient $\thestack^{ss}(l)/G$ resulting from choosing $l \in M^{\mathcal{W}}_{\mathbb{R}}$ positive or negative respectively under the identification with $\mathbb{R}$.
Recall that the window subcategories induce \textit{window-shift} autoequivalences by composing along the following set of equivalences.
\begin{equation*}
\begin{tikzcd}
	& {\mathscr{W}_{j}} && {\mathscr{W}_{k}} \\
	{\db(X_+)} && {\db(X_-)} && {\db(X_+)}
	\arrow[from=1-2, to=2-3]
	\arrow[from=1-4, to=2-5]
	\arrow[from=2-1, to=1-2]
	\arrow[from=2-3, to=1-4]
\end{tikzcd}
\end{equation*}
Here we have made a choice of two distinct window subcategories, and used the fact that each window is derived equivalent to both GIT quotients.

We are interested in describing the autoequivalences of $\db(X_+)$ as the twist around a spherical functor in the case where $k = j + 1$, i.e. the windows come from adjacent one-dimensional faces $C_j$ and $C_{j + 1}$.
The first result towards this is the following description of the window-shift equivalence $\mathscr{W}_{j} \to \db(X_-) \to \mathscr{W}_{j + 1}$.

\begin{thm}[{\cite[Theorem 5.12]{schobers}}]\label{thm:sod}
    There is an equivalence $\mathscr{W}_{j} \to \mathscr{W}_{j + 1}$ given by mutation in a semiorthogonal decomposition of $\mathscr{E}_{D_j}$:
    \[
    \mathscr{E}_{D_j} = \langle \mathscr{K}_{j}, \mathscr{W}_{j}\rangle = \langle \mathscr{W}_{j + 1}, \mathscr{K}_{j}\rangle
    \]
    over a subcategory $\mathscr{K}_{j} \subset \db(\thestack)$.
\end{thm}
The mutation equivalence can be described as including the window into the larger subcategory and projecting across the orthogonal $\mathscr{K}_j$.

\subsubsection{Orthogonal subcategories}
The orthogonal subcategories $\mathscr{K}_{j}$ appearing in \Cref{thm:sod} are generated by certain objects supported on the unstable locus for $X_{-}$.
Let $\lambda$ be a one-parameter subgroup lying in $N^{-}_{\mathbb{R}}$.
Denote by $Z_{\lambda}$ the fixed locus of $\lambda$ and $Y_{\lambda} \subset \thestack$ the subspace 
\[
Y_{\lambda} := \{ x \in \thestack \;\vert\; \lim_{t \to 0} \lambda(t) \cdot x \in Z_{\lambda} \}.
\]
This is the locus of points destabilized by $\lambda$, and the space $S_{\lambda} := GY_{\lambda}$ is the locus destabilized by the conjugates of $\lambda$.

Let $\widetilde{S}_{\lambda} := G \times_{B} Y_{\lambda} \cong Y_{\lambda}/B$; this is a subspace of $\thestack \times G/B$.
We will take $G = GL(V)$ and $B$ to be the subgroup of upper triangular matrices, and henceforth identify $G/B \cong \mathbb{P}(V)$.
We let
\[
\sigma : \thestack \times \mathbb{P}(V) \to \thestack
\]
be the projection.
We have the following `Springer-type resolution' diagram.
\[\begin{tikzcd}
	{\widetilde{S}_{\lambda}} && {\thestack \times \mathbb{P}(V)} && {\mathbb{P}(V)} \\
	\\
	{S_{\lambda}} && \thestack && pt
	\arrow[hook, from=1-1, to=1-3]
	\arrow[from=1-1, to=3-1]
	\arrow[from=1-3, to=1-5]
	\arrow["\sigma"', from=1-3, to=3-3]
	\arrow[from=1-5, to=3-5]
	\arrow[hook, from=3-1, to=3-3]
	\arrow[from=3-3, to=3-5]
\end{tikzcd}\]
The space $\widetilde{S}_{\lambda}$ is a $G$-equivariant vector bundle over $\mathbb{P}(V)$.

\begin{figure}[!h]
    \centering
    \includegraphics[width=0.5\linewidth]{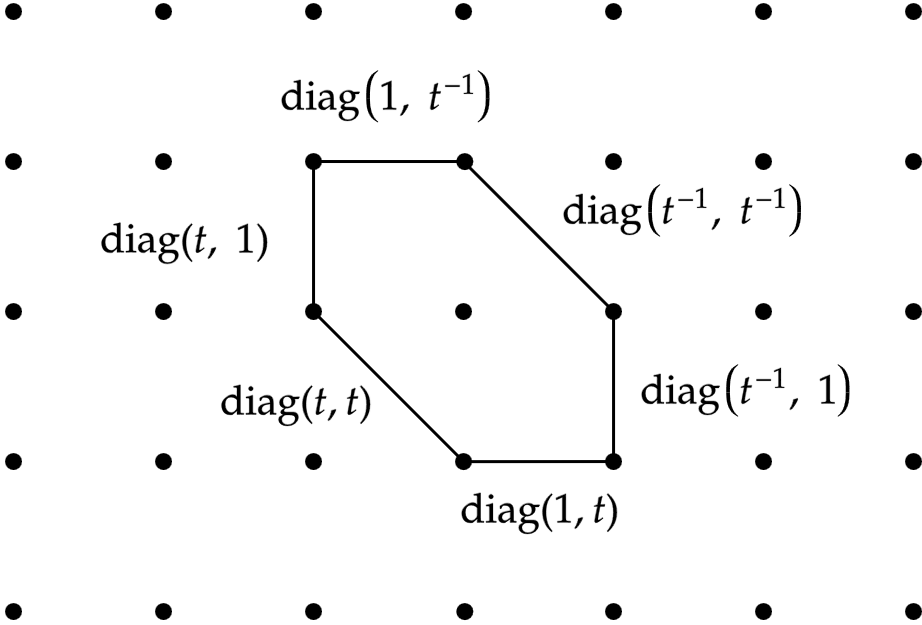}
    \caption{One-parameter subgroups corresponding to each face of $\overline{\nabla}$.}
    \label{fig:1ps}
\end{figure}

Next, consider a facet $F$ of $\overline{\nabla}_{D_j}$ and let $\mu_{F} \subset N_{\mathbb{R}}$ be the set of $\lambda$ such that $\langle \lambda, f \rangle \geq u$ for some $u$ and all $f \in \overline{\nabla}_{D_{j}}$, with equality if and only if $f \in F$.
In general, the effect of a window-shift may not be concentrated along all unstable strata; the purpose of the set $\mu_{F}$ is to determine precisely which strata are relevant for the particular window-shift.
For our polytope, the sets $\mu_{F}$ are rays for a given one-dimensional face $F$, and we display the one-parameter subgroups generating the ray in \Cref{fig:1ps}.

\begin{thm}[{\cite{hlsam}, \cite[Corollary A.3]{schobers}}]\label{thm:orthogonal}

    Take vectors $\xi \in C_{j}$ and $\xi^\prime \in D_j$ and let $\epsilon = \xi - \xi^{\prime}$.
    Let $\nu_{\epsilon} := \{\lambda \in N^{-}_{\mathbb{R}} \; \vert \; \langle \lambda, \epsilon\rangle > 0\}$.
    The subcategory $\mathscr{K}_{j}$ is generated by the objects
    \begin{align*}
        \{ \sigma_{*}(\sheaf_{\widetilde{S}_{\lambda}} \otimes \chi) \; \vert \; \chi \in \overline{\nabla}_{D_j} \setminus \overline{\nabla}_{C_{j}}, \; \lambda \in \nu_{\epsilon} \cap \mu_{F} \; \text{for all faces} \; F \; \text{containing} \; \chi\}
    \end{align*}
\end{thm}
Note that in this theorem, we consider $\chi$ as a weight of $B$; since $\widetilde{S}_{\lambda}$ is a bundle over $G/B$, the notation $\sheaf_{\widetilde{S}_{\lambda}} \otimes \chi$ refers to taking the line bundle on $G/B$ associated to $\chi$ and pulling back to $\widetilde{S}_{\lambda}$.
We highlight that the purpose of requiring $\lambda$ to pair positively with $\epsilon$ is to ensure that only the unstable stratum for $X_{-} = \thestack^{ss}(\epsilon)/G$ are considered.

The autoequivalences resulting from mutation can then be written as spherical twists using the following.

\begin{thm}[{\cite[Lemma 6.7]{hlsam},\cite[Proposition 3.4]{hlshipman}}]
    The window-shift autoequivalence 
    \[\db(X_+) \to \mathscr{W}_{j} \to \mathscr{W}_{j+1} \to \db(X_+)\] 
    is the twist around the composition
    \[
    \mathscr{K}_{j} \to \mathscr{E}_{D_j} \to \db(X_+)
    \]
    where the latter map is the restriction functor.
\end{thm}

\subsubsection{Unstable strata} \label{subsubsec:strata}
Here we analyse the unstable strata for $X_-$.
On $\mathbb{P}(V)$, we denote by $L$ the tautological subspace bundle, and $Q := V/L$ the quotient bundle.
Furthermore, we use $D := \operatorname{det} V$ to avoid confusion with the $\sheaf(1)$ line bundle on $\mathbb{P}(V) \cong \mathbb{P}^1$.
We will also reuse these labels for the pullbacks of these bundles to $\thestack \times \mathbb{P}(V)$.

There is a two-level Kempf-Ness stratification, with the deepest stratum induced by the one-parameter subgroup
\[
\lambda_0(t) = \operatorname{diag}(t^{-1}, t^{-1})
\]
given by 
\[
S_0 := \{\alpha = 0\}.
\]
After removing this locus, we consider the one-parameter subgroup
\[
\lambda_1(t) = \operatorname{diag}(1, t^{-1})
\]
which destabilizes the locus where both $\beta$ and $\gamma$ share $\alpha$ as an eigenvector.
We denote this by $S^{o}_1 \subset \thestack \setminus S_0$, and its closure in $\thestack$ by $S_1$.
The subscheme $S_1$ is singular, and $\widetilde{S}_1 := Y_{\lambda_1} \times_{B} G \subset \thestack \times \mathbb{P}(V)$ gives a resolution. 
This can be described as choosing a line $L \subset V$ as the common eigenvector of $\beta$ and $\gamma$, and enforcing $\alpha \subset L$.
Then, $\widetilde{S}_1$ is the following vector bundle over $\mathbb{P}(V)$:
\[
\widetilde{S}_1 \cong \operatorname{Tot}\left(L_p \oplus V^*_{\alpha^*} \oplus \operatorname{Hom}(L \oplus Q, L)^{\oplus 2}_{\beta, \gamma} \right)_{\mathbb{P}(V)}.
\]
Similarly, we denote by $\widetilde{S}_0$ the locus $\{\alpha = 0\} \cong S_0 \times \mathbb{P}(V) \subset \thestack \times \mathbb{P}(V)$.

\subsubsection{Orthogonal subcategories for the universal flop} The Picard group acts as tensoring by $\sheaf(1)$, and so identifies $\mathscr{W}_{j}$ with $\mathscr{W}_{j + 2}$.
Working up to this action, we only need to consider the two window shifts $\mathscr{W}_{-2} \to \mathscr{W}_{-1}$ and $\mathscr{W}_{-1} \to \mathscr{W}_{0}$ in order to understand the monodromy action.
Let $\overline{\nabla}_{\text{fib}}$ and $\mathscr{E}_{\text{fib}}$ be the polytope and categories associated to the zero dimensional face $D_{-2} = \{-\tfrac{1}{2}\}$, and $\overline{\nabla}_{\text{con}}$ and $\mathscr{E}_{\text{con}}$ be the polytope and categories associated to the face $D_{-1} = \{0\}$.

By \Cref{thm:sod}, there are then semiorthogonal decompositions
\begin{align*}
    \mathscr{E}_{\text{fib}} & = \langle \mathscr{K}_{\text{fib}}, \mathscr{W}_{-2} \rangle = \langle \mathscr{W}_{-1}, \mathscr{K}_{\text{fib}} \rangle \\
    \mathscr{E}_{\text{con}} & = \langle \mathscr{K}_{\text{con}}, \mathscr{W}_{-1} \rangle = \langle \mathscr{W}_{0}, \mathscr{K}_{\text{con}} \rangle.
\end{align*}

\begin{prop}
    The category $\mathscr{K}_{\text{fib}}$ is generated by $\sheaf_{S_0}$.
\end{prop}
\begin{proof}
    The only character in $\overline{\nabla}_{D_{-2}} \setminus \overline{\nabla}_{C_{-1}}$ corresponds to $\sheaf$, and by considering the face of $\overline{\nabla}_{D_{-1}}$ containing the character, we see that the left orthogonal is then generated by the pushforward of the structure sheaf on $S_0$.
\end{proof}

\begin{prop}\label{prop:Kcon}
    The category $\mathscr{K}_{\text{con}}$ is generated by the objects $\mathcal{F} := \sigma_{*}\sheaf_{\widetilde{S}_1}(Q)$ and $\mathcal{G} := \sheaf_{S_0}(V)$.
\end{prop}
\begin{proof}
    The characters in $\overline{\nabla}_{D_{-1}} \setminus \overline{\nabla}_{C_0}$ correspond to the bundles $V, \Sym^2V(-1)$.
    The character for $V$ lies in the closure of two codimension one faces of $\overline{\nabla}_{\text{con}}$, corresponding to the one-parameter subgroups $\lambda_0(t) = \operatorname{diag}(t^{-1}, t^{-1})$ and $\lambda_1 = \operatorname{diag}(1, t^{-1})$.
    Similarly, $\Sym^2V(-1)$ lies in the closure of faces corresponding to the one-parameter subgroups $\lambda_1 = \operatorname{diag}(1, t^{-1})$ and additionally $\lambda^{\prime} = \operatorname{diag}(t, 1)$.
    We discard the latter since it pairs negatively with a vector $\epsilon \in C_{-1}$.

    Hence, the left orthogonal subcategory to $\mathscr{W}_{-1}$ in $\mathscr{E}_{\text{con}}$ is generated by the objects 
    \[\mathscr{K}_{\text{con}} := \langle \sheaf_{S_0}(V), \sigma_{*}\sheaf_{\widetilde{S}_1}(Q), \sigma_{*}\sheaf_{\widetilde{S}_1}(Q^2 \otimes D)\rangle.\]
    However, since $\mathscr{E}_{\text{con}}$ is generated by the bundles $\sheaf, V, V^*$ and $ \Sym^2V(-1)$, by considering the locally free resolutions of $\mathcal{F}$ and $\mathcal{G}$ as computed in \Cref{lemma:resF,lemma:resG}, we see that it is sufficient to take $\mathscr{W}_{0} = \langle \sheaf, V \rangle$ and $\mathcal{F}, \mathcal{G}$ as generators for $\mathscr{E}_{\text{con}}$, and consequently $\mathscr{K}_{\text{con}}$ can be generated by $\mathcal{F}, \mathcal{G}$.
\end{proof}

\begin{remark}
    Another way to see that we do not need to take all generators is to consider the following.
    Take the sheaf given by $\mathcal{G}/[\beta, \gamma] := \operatorname{coker} [\beta, \gamma] : \mathcal{G} \to \mathcal{G}$. 
    This sheaf is supported on $S_0 \cap S_1$, and specifically we see that it has rank $1$ generically and rank $2$ where $\beta, \gamma$ share all eigenvectors.
    Then we note that on $\widetilde{S}_1$, the cone on the map $\sheaf \to L \cong Q^{-1} \otimes D$ cuts out the intersection with $\widetilde{S}_0$.
    Hence, there is a map $\sigma_{*}\sheaf_{\tilde{S}_1}(Q^{2} \otimes D^{-1}) \to \sigma_{*}\sheaf_{\tilde{S}_1}(Q)$ with cone precisely given by $\mathcal{G}/[\beta, \gamma]$.
\end{remark}

The main objective now is to provide a description of the categories $\mathscr{K}_{\text{fib}}$ and $\mathscr{K}_{\text{con}}$, and relate the resulting spherical functors to the Donovan--Wemyss fibre and contraction algebra twists.

\subsection{The fibre twist}\label{subsec:fibre}

As suggested by our choice of notation, the claim of this section is that the spherical functor with source $\mathscr{K}_{\text{fib}}$ gives an autoequivalence which agrees with the twist around the functor $\db(\Afib) \to \db(X_+)$.

\begin{lemma}\label{lemma:endOs0}
    We have that $\End(\sheaf_{S_0}) \cong \Hom(\sheaf_{S_0}, \sheaf_{S_0}) \cong \mathbb{C}[\tr{\beta^2}, \tr \gamma^2, \tr \beta \gamma]$ is a polynomial ring in three variables.
\end{lemma}
\begin{proof}
    Firstly, we have that $\Gamma(S_0, \sheaf) \cong \CC[\tr \beta^2, \tr \gamma^2, \tr \beta\gamma]$. 
    This follows from noting that $S_0$ is the stack $[V^* \oplus \Sym^2V(-1)^{\oplus 2} / GL(V)]$ and using a classical result for the invariant ring, see for example \cite[p. 21]{kraft-procesi}.
    The fact that $\End_{\thestack}(\sheaf_{S_0}) \cong \End_{S_{0}}(\sheaf)$ follows from \cite[Amplification 3.27]{hl} since the conormal bundle of $S_0$ is given by $V^*$, which has strictly positive $\lambda_0$-weights.
\end{proof}

Now note that the deformation parameters $T_0^{\beta}, T_0^{\gamma}, T_0^{\delta}$ in the original quiver correspond to elements of $\Gamma(\thestack, \sheaf)$.
Then, as $\Lambda_{-1}$ modules, it is clear that there is an identification $\Gamma(S_0, \sheaf) \cong \Afib$, since one can verify that the deformation parameters act as $T_0^{\beta} = \tr \beta^2, T_0^{\gamma} = \tr \gamma^2, T_0^{\delta} = \tr \beta^2 + 2 \tr \beta\gamma + \tr \gamma^2$ up to scaling.

\begin{prop}\label{prop:afibfunctor}
    There is a window-shift autoequivalence $\omega_{\text{fib}}$ given by the twist around a spherical functor
    \[\db(\mathbb{C}[\tr{\beta^2}, \tr \gamma^2, \tr \beta \gamma]) \to \db(X_+).\]
    Moreover, using the identification $\mathbb{C}[\tr{\beta^2}, \tr \gamma^2, \tr \beta \gamma] \cong \Afib$, this recovers the functor given by the composition
    \[
    \db(\Afib) \to \db(\Lambda_{-1}) \to \db(X_+).
    \]
\end{prop}
\begin{proof}
    The first statement is just a corollary of \Cref{lemma:endOs0}.
    The restriction of $\sheaf_{S_0}$ to $X_+$ is $\sheaf_{E_1}$, and $\Rhom_{X_+}(V^* \oplus \sheaf, \sheaf_{E_1}) \cong \Rhom_{E_1}(V^* \oplus \sheaf, \sheaf).$
    To compute this, note that $E_1$ is the GIT quotient of the stack $S_0 \cong [V^* \oplus \Sym^2V(-1)^{\oplus 2} / GL(V)]$.
    The bundles $\sheaf, V^*$ then satisfy appropriate grade restriction rules by essentially the same calculation as for the full stack $\thestack$, so applying \cite[Theorem 2.10]{hl} we find that it suffices to compute morphisms on the stack.
    But here we know that $\Rhom_{S_0}(\sheaf, \sheaf) \cong \Afib$, and we also have that 
    \begin{align*}
        \Rhom_{S_0}(V^*, \sheaf) \cong Hom_{GL(V)}(V^*, \Sym^{\bullet}(V \oplus \Sym^2V(-1)^{\oplus 2})) = 0.
    \end{align*}
\end{proof}

Letting $\operatorname{Spec} R$ be the singular base of the universal flop, recall that there is an irreducible component of the singular locus $Z_1$, and the coordinate ring is identified with $\Afib$ as an $R$-module.
We can then recover the geometric picture.
\begin{prop}
    The window-shift autoequivalence $\omega_{\text{fib}}$ is naturally isomorphic to the twist around the functor given by $j_*\pi^*$, defined using the following diagram.

    \begin{equation*}
    \begin{tikzcd}
    	& {E_1} \\
    	{\mathbb{A}^3 \cong Z_1} && {X_+}
    	\arrow["\pi"', from=1-2, to=2-1]
    	\arrow["j", hook, from=1-2, to=2-3]
    \end{tikzcd}
    \end{equation*}

\end{prop}
\begin{proof}
    We give a direct proof of this in our setting, similarly to \cite[Proposition 2.2]{ds14}, but a more abstract proof along the lines of \cite[Proposition 3.4]{hlshipman} is also possible.
    Consider the following diagram.
    \[\begin{tikzcd}
	{S_0} && \thestack \\
	& {E_1} && {X_+} \\
	{\mathbb{A}^3} & {\mathbb{A}^3}
	\arrow["{\hat{\jmath}}", hook, from=1-1, to=1-3]
	\arrow["{\hat{\pi}}", from=1-1, to=3-1]
	\arrow["r", hook', from=2-2, to=1-1]
	\arrow["j", hook, from=2-2, to=2-4]
	\arrow["\pi"', from=2-2, to=3-2]
	\arrow["r", hook', from=2-4, to=1-3]
	\arrow[equals, from=3-1, to=3-2]
\end{tikzcd}\]
    Let $f := j_*\pi^{*}$ and $\hat{f} := \hat{\jmath}_*\hat{\pi}^*$.
    Let $r_{j} : \mathscr{W}_{j} \to \db(X_+)$ be the restriction equivalence.
    To perform the window-shift autoequivalence, we take the cone on the map $r_{-1} \circ \hat{f} \circ \hat{f}^{R} \circ r^{-1}_{-2} \to Id$.
    But note that for an object in $\mathscr{W}_{-2}$, we have that $\hat{f}^{R} \vert_{\mathscr{W}_{-2}} \simeq f^{R} \circ r_{-2}$ 
    To check this, we note that the vector bundles $\sheaf, V^*$ on $S_0$ have no higher cohomology upon restricting to $E_1$.
    Similarly, $r_{-1} \circ \hat{f} \vert_{\mathscr{W}_{-1}} \simeq f$, where now we use that since $E_1$ is the intersection of $S_0$ with the semistable locus for $X_+$, the top right square is a pullback diagram.

\end{proof}

\begin{remark}
    In fact, it is possible to fully build the autoequivalence geometrically, using methods closer to those in \cite{ds14}.
    Let $\pi : E_1 \to Z_1$ and $j : E_1 \to X_+$.
    The functor $F = j_*\pi^*$ can then be shown to be spherical; the key is to show that the co-twist $F^LF$ is an autoequivalence of $\db(\mathbb{A}_3)$.
    Generically along the $E_1$ component, the autoequivalence is a family of Atiyah flops, and in fact one can show that the resulting cotwist is also a shift in $\db(\mathbb{A}_3)$, as for the Atiyah flop, with the proof also being similar.
\end{remark}

\subsection{The contraction twist}\label{sec:contractiontwist}

The main claim of this section is that the autoequivalence given by the twist around the functor provided by the machinery of \Cref{sec:algorithm}
\[
    \mathscr{K}_{\text{con}} \to \db(X_+)
\]
is naturally isomorphic to the autoequivalence given by the twist around Donovan--Wemyss's spherical functor $\db(\Acon) \to \db(X_+)$.
Recall the following presentation of $\Acon$:
\begin{align*}
    \Acon \cong \frac{\mathbb{C}[t]\langle\beta, \gamma\rangle}{\begin{matrix}
        [\beta^2, \gamma], [\gamma^2, \beta], t [\beta, \gamma]
    \end{matrix}}.
\end{align*}
Loosely speaking, we may think of this as the coordinate ring of a non-commutative scheme with two irreducible components, one of which has coordinate ring
\[
\Acon / t \cong \frac{\CC\langle \beta, \gamma \rangle}{\begin{matrix}
    [\beta^2, \gamma],
    [\gamma^2, \beta]
\end{matrix}}
\]
and the other with coordinate ring
\[
\Acon / [\beta, \gamma] \cong \CC[t, b, c].
\]
Note here that $\beta\gamma - \gamma\beta$ anti-commutes with $\beta, \gamma$, so defines a two-sided ideal.
During the proof of \Cref{lemma:Phi} we will show that the modules $\Acon/t$ and $\Acon/[\beta, \gamma]$ have infinite projective dimension, hence $\Acon$ is singular.

Recall that by \Cref{prop:Kcon}, $\mathscr{K}_{\text{con}}$ is generated by two objects $\mathcal{F}$ and $\mathcal{G}$.
We therefore compute the endomorphism algebra of $\mathcal{F} \oplus \mathcal{G}$ on $\thestack$.

\begin{lemma}\label{lemma:endf}
We have that
\[
\operatorname{End}(\mathcal{F}) \cong \Hom(\mathcal{F}, \mathcal{F}) \cong \mathbb{C}[t, b, c].
\]
\end{lemma}
\begin{proof}
    Consider the stack $\mathfrak{V} := \thestack \setminus S_0$.
    Denote by $\hat{\mathcal{F}}$ the restriction of $\mathcal{F}$ to this locus.
    Note that the sheaf $\hat{\mathcal{F}}$ is the pushforward of the structure sheaf of the locus $S^{o}_1 \subset \mathfrak{V}$, which is a vector bundle over the fixed locus of the one-parameter subgroup $\lambda_1(t) := \operatorname{diag}(1, t^{-1})$ in $\mathfrak{V}$.
    This is given by the stack 
    \[\{\beta_{01} = \gamma_{01} = \beta_{10} = \gamma_{10} = \alpha_1 = \alpha_1^* = 0, \alpha \neq 0\} = [\Aff^4 / \CC^*] \cap \{\alpha \neq 0\}\] and is isomorphic to $\mathbb{A}^3$.
    The coordinates on this space are $t := \alpha_0\alpha_0^*, b := \alpha_0\beta_{00}, c := \alpha_0\gamma_{00}$, where we recall that we label $t = \alpha\alpha^* \in \Gamma(\thestack, \sheaf)$.
    By \cite[Corollary 3.28]{hl} we have that $\End(\hat{\mathcal{F}}) \cong \Gamma(S^{o}_1, \sheaf) \cong \CC[t, b, c]$.

    Next, we take the locally free resolution of $\mathcal{F}$ as given in \Cref{lemma:resF}.
    Restricting this resolution to the fixed locus of $\lambda_0(t)$, we obtain 
    \[
        \sheaf^{\oplus 3}
        \to \sheaf^{\oplus 2} \oplus \sheaf(-1)^{\oplus 2} \oplus \sheaf
          \to \sheaf(-1)^{\oplus 2}
    \]
    where we are temporarily using $\sheaf(-1)$ to denote the line bundle of weight $-1$ under $\lambda_0$.
    In particular, the $\lambda_0$ weights lie in the range $[-1, 0]$, so by derived Kirwan surjectivity \cite[Theorem 2.10]{hl} we find that $\End_{\thestack}(\mathcal{F}) \cong \End_{\mathfrak{V}}(\hat{\mathcal{F}})$.
\end{proof}

We emphasise that the polynomial ring $\mathbb{C}[t, b, c]$ is distinct, e.g. as an $R$-module, to the polynomial ring in three variables $\Afib$ that appeared in the previous section.

\begin{lemma}\label{lemma:endg}
We have that
\[
\operatorname{End}(\mathcal{G}) \cong \Hom(\mathcal{G}, \mathcal{G}) \cong \mathbb{C}\langle \beta, \gamma\rangle/([\beta^2, \gamma], [\gamma^2, \beta]).
\]
\end{lemma}
\begin{proof}
    Using \cite[Amplification 3.27, Proposition 3.25]{hl}, since the conormal bundle $N^{\vee}_{S_0/\thestack}$ restricted to $\{\alpha = \alpha^* =0\}$ has positive weights we find that $\End_{\thestack}(\mathcal{G}) \cong \End_{S_0}(V)$.
    Recall that 
    \[S_0 \cong [V^* \oplus (\Sym^2V(-1))^{\oplus 2} / GL(V)].
    \]
    To compute $\End(V)$ on this stack, it is equivalent to compute $\End(V)$ on the stack \[[\Sym^2V(-1)^{\oplus 2} / GL(V)].\]
    This is now a classical result in the theory of trace rings; we can use \cite[Theorem I.3.1]{bruyntrace} combined with \cite[Theorem III.5.1]{bruyntrace}, to see that the desired algebra is the Ore extension
    \[
        \mathbb{C}[a_{11}, a_{12}, a_{22}][x_1][x_2, \sigma_2, \delta_2]
    \]
    where $\sigma_{2}(x_1) = -x_1$ and $\delta_{2}(x_1) = 2a_{12}$ and the action on the other elements is trivial.
    The resulting $\mathbb{C}[a_{11}, a_{12}, a_{22}]$-algebra is therefore generated by $x_1, x_2$  with the relations
    \begin{align*}
        x_2 x_1 + x_1 x_2 = 2a_{12} && x_1^2 = a_{11} && x_2^2 = a_{22}
    \end{align*}
    which is isomorphic under $x_1 \mapsto \beta, x_2 \mapsto \gamma$ to our desired algebra.
\end{proof}

The next lemma shows that the objects $\mathcal{F}, \mathcal{G}$ are in fact semiorthogonal.

\begin{lemma}
We have that
\[
\operatorname{Ext}^{i}(\mathcal{G}, \mathcal{F}) = 0
\]
for all $i \in \mathbb{Z}$.
\end{lemma}
\begin{proof}
    Pulling back the locally free resolution of $\mathcal{G}$ to $\widetilde{S}_1$ and then applying $\underline{\operatorname{Hom}}(-, Q)$ we obtain:
    \[
    0 \to V^* \otimes Q \to (\sheaf \oplus \Sym^2V(-1)) \otimes Q \to V \otimes Q \to 0.
    \]
    On $\widetilde{S}_1$, this complex has cohomology at the middle term and at the final term.

    At the final term, we have to compute global sections of $V \otimes Q$ on the space $\widetilde{S}_0 \cap \widetilde{S}_1$.
    This is a vector bundle over $\mathbb{P}^1$, and pushing down to $\mathbb{P}^1$, we see that we need to count copies of $V^*$ in the $GL(V)$ representation
    \[
    \Gamma\left(\mathbb{P}^1, Q \otimes \Sym^{\bullet}( \sheaf^{\oplus 2} \oplus (Q^2 \otimes D^{-1})^{\oplus 2} \oplus V)\right).
    \]
    Hence there is no cohomology in degree $2$.
    On $\widetilde{S}_1$, the kernel of the map $\alpha : V^* \to \sheaf$ is given by $Q^{-1} \to V^*$, so the kernel at the middle term is $Q^{-1} \otimes V \otimes Q \cong V$.
    Hence we need to compute global sections of $V$ on the intersection $\widetilde{S}_0 \cap \widetilde{S}_1$.
    So now we conclude by observing that
    \begin{align*}
        & \Gamma(\widetilde{S}_0 \cap \widetilde{S}_1, V) \\
        \cong & \Hom_{GL(V)}(V^*, \Gamma(\mathbb{P}(V), \Sym^{\bullet}(\sheaf^{\oplus 2} \oplus (Q^2 \otimes D^{-1})^{\oplus 2} \oplus V))) = 0.
    \end{align*}
\end{proof}

\begin{lemma}\label{lemma:bimodule}
We have that
\[
\operatorname{Ext}^{\bullet}(\mathcal{F}, \mathcal{G}) \cong 
\operatorname{Ext}^1(\mathcal{F}, \mathcal{G}) \cong 
      \mathbb{C}[b, c]
\]
and furthermore, the algebra is isomorphic as a right module over $\operatorname{End}(\mathcal{F})$ to the algebra $\operatorname{End}(\mathcal{F})/t$, while as a left module over $\operatorname{End}(\mathcal{G})$ it is isomorphic to $\operatorname{End}(\mathcal{G})/[\beta, \gamma]$.
\end{lemma}
\begin{proof}
    We begin by computing $\sigma^!\mathcal{G}$, where $\sigma : \thestack \times \mathbb{P}^1 \to \thestack$.
    The relative canonical bundle is given by
    \[
    \omega_\sigma \cong Q^{-1} \otimes L \cong Q^{-2} \otimes D.
    \]
    Hence we want to compute
    \[
    \operatorname{Ext}^{\bullet}(\sheaf_{\widetilde{S}_1}(Q), \sheaf_{\widetilde{S}_0}(V \otimes Q^{-2} \otimes D)[1]) \cong \operatorname{Ext}^{\bullet}(\sheaf_{\widetilde{S}_1}, \sheaf_{\widetilde{S}_0}(V \otimes Q^{-3} \otimes D)[1]).
    \]
    Now we use the locally free resolution for $\sheaf_{\widetilde{S}_1}$ obtained in \Cref{lemma:resF}.
    \[
    0 \to (Q^{-5} \otimes D^2) \to (Q^{-3} \otimes D) \oplus (Q^{-4} \otimes D^2) \to (Q^{-2} \otimes D)^{\oplus 2} \oplus Q \to \sheaf \to 0.
    \]
    Ignoring the twist by $V \otimes Q^{-2} \otimes D$ momentarily, we apply $\underline{\Hom}(-, \sheaf_{\widetilde{S}_0})$ to this sequence to obtain
    \[
    0 \to \sheaf \to (Q^2 \otimes D^{-1})^{\oplus 2} \oplus Q \to (Q^{3} \otimes D^{-1}) \oplus (Q^{4} \otimes D^{-2}) \to  Q^{5} \otimes D^{-2} \to 0.
    \]
    We can present the maps in this sequence as follows.
    We use the notation of \Cref{lemma:resF}, where the maps $s_{\beta}, s_{\gamma} \in \Hom(L, Q)$ are induced by the compositions $L \to V \to V \to Q$.
    The first map is 
    \[
    \xymatrix@C=6pc{
        \sheaf 
        \ar[r]^-{
            \begin{bsmallmatrix} 
                s_\beta & s_\gamma & 0
            \end{bsmallmatrix}
        } & \left( Q^{2} \otimes D^{-1} \right)^{\oplus 2} \oplus Q
    }
    \]
    The second map is
    \begin{equation*}
    \xymatrix@C=8pc{
      (Q^2 \otimes D^{-1})^{\oplus 2} \oplus Q \ar[r]^-{\begin{bsmallmatrix}
        0 & 0 & s_{\gamma}  \\ 0 & 0 & -s_{\beta} \\
      -s_{\beta} & -s_{\gamma} & 0\end{bsmallmatrix}} &
      (Q^3 \otimes D^{-1})^{\oplus 2} \oplus (Q^4 \otimes D^{-2})
    }.
    \end{equation*}
    The final map is
    \begin{equation*}
    \xymatrix@C=6pc{
      (Q^3 \otimes D^{-1})^{\oplus 2} \oplus (Q^4 \otimes D^{-2}) \ar[r]^-{\begin{bsmallmatrix}
          -s_{\gamma} \\ s_\beta \\ 0
      \end{bsmallmatrix} } &
      Q^5 \otimes D^{-2}
    }.
    \end{equation*}
    In degree $3$, the cohomology is $Q^5 \otimes D^{-2}$ on the intersection $\widetilde{S}_0 \cap \widetilde{S}_1$, since the sections $s_{\gamma}, s_{\beta}$ cut out the locus $\widetilde{S}_0 \cap \widetilde{S}_1$ in $\widetilde{S}_0$.
    Remembering that we still have to tensor by $V \otimes Q^{-3} \otimes D$, we see that in degree $3$ we need global sections of $V \otimes Q^2 \otimes D^{-1}$ on $\widetilde{S}_1 \cap \widetilde{S}_0$, but one observes that:
    \begin{align*}
        & \; \Gamma(\widetilde{S}_0 \cap \widetilde{S}_1, V \otimes Q^2 \otimes D^{-1}) \\ 
        \cong & \; \Hom_{GL(V)}\left(V^*, \Gamma(\mathbb{P}^1, Q^2 \otimes D^{-1} \otimes \Sym^{\bullet}(\sheaf^{\oplus 2} \oplus (Q^2 \otimes D^{-1})^{\oplus 2} \oplus V))\right) = 0.
    \end{align*}
    In degree $2$ we need to compute global sections of $V \otimes (Q \otimes D^{-1})$ on $\widetilde{S}_0 \cap \widetilde{S}_1$.
    This gives a polynomial ring $\mathbb{C}[b, c]$ by taking the symmetric algebra on the coordinates of the trivial line bundles:
    \begin{align*}
        & \; \Gamma(\widetilde{S}_0 \cap \widetilde{S}_1, V \otimes Q \otimes D^{-1}) \\ 
        \cong & \; \Hom_{GL(V)}(V^*, \Gamma(\mathbb{P}^1, Q \otimes D^{-1} \otimes \Sym^{\bullet}(\sheaf^{\oplus 2} \oplus (Q^2 \otimes D^{-1})^{\oplus 2} \oplus V)))
        \\ \cong & \; \CC[b, c].
    \end{align*}
    Recalling the shift of $1$ induced by $\sigma^!$, we therefore have the first statement of the lemma.

    Using the fact that $\widetilde{S}_1 \to S_1$ is an isomorphism outside of codimension $2$, we can identify the elements $b, c \in \End(\mathcal{F})$ with sections in $\Gamma(\widetilde{S}_1, \sheaf)$, coming from the trivial summands of 
    \[\widetilde{S}_1 \cong \operatorname{Tot}\left( \sheaf^{\oplus 2}_{b, c} \oplus (Q^{-2} \otimes D)^{\oplus 2}_{b', c'} \oplus (Q^{-1} \otimes D)_{p} \oplus V^*_{\alpha^*} \right)_{\mathbb{P}(V)}\]
    while $t \in \End(\mathcal{F})$ can be identified with the section $\alpha^*p \in \Gamma(\widetilde{S}_1, \sheaf)$.
    The module structure of $\operatorname{Ext}^{1}(\mathcal{F}, \mathcal{G})$ over $\operatorname{End}(\mathcal{F})$ then follows directly from the module structure of \[\Gamma(\widetilde{S}_1, \sheaf_{\widetilde{S}_0 \cap \widetilde{S}_1}(V \otimes Q \otimes D^{-1}))\] over $\Gamma(\widetilde{S}_1, \sheaf)$. 
    Similarly, the module structure over $\operatorname{End}(\mathcal{G})$ follows by considering 
    \[\Gamma(\widetilde{S}_0, \sheaf_{\widetilde{S}_0 \cap \widetilde{S}_1}(V \otimes Q \otimes D^{-1}))\] as a module over $\End_{\widetilde{S}_0}(V) \cong \End(\mathcal{G})$.
\end{proof}

The subcategory of $\db(\thestack)$ generated by $\mathcal{F}$ and $\mathcal{G}$ acts as a source category for our spherical functor; the next proposition tells us that this category is smooth.

\begin{prop}\label{prop:smooth-algebra}
The algebra $\End(\mathcal{F} \oplus \mathcal{G}[1])$ has finite global dimension.
\end{prop}
\begin{proof}
    The algebra $\End(\mathcal{F})$ is a polynomial ring and the fact that $\End(\mathcal{G})$ has finite global dimension is a classical result, see for example \cite[Theorem II.3.2]{bruyntrace}. 
    It is then straightforward to see, using the description in \Cref{lemma:bimodule}, that the bimodule $\Ext^{1}(\mathcal{F}, \mathcal{G})$ is perfect over $\End(\mathcal{F}) \otimes \End(\mathcal{G})^{\operatorname{op}}$, so taking into account the semiorthogonality and viewing $\End(\mathcal{F} \oplus \mathcal{G}[1])$ as the upper triangular algebra given by the gluing of two algebras along a bimodule, we may conclude by \cite[Proposition 3.11]{lunts}.
\end{proof}

At this point, it remains to be shown that the resulting autoequivalence is the contraction algebra twist given by Donovan--Wemyss.
Firstly, we show that the functor constructed admits a map to $\db(\Acon)$. 

\begin{lemma}
    Let $\mathcal{E} := \operatorname{Cone}(\mathcal{F} \to \mathcal{G}[1])[-1]$, where we take the cone on the map corresponding to $1 \in \operatorname{Ext}^1(\mathcal{F}, \mathcal{G}) \cong \mathbb{C}[b, c]$.
    Then $\End(\mathcal{E}) \cong \Hom(\mathcal{E}, \mathcal{E}) \cong \Acon$, and hence there is a functor $\Phi : \db(\mathcal{F} \oplus \mathcal{G}) \to \db(\Acon)$ given by
    \[
    \Phi = \Rhom(\mathcal{E}, -).
    \]
\end{lemma}
\begin{proof}
    The endomorphisms of $\mathcal{E}$ are given by the fiber product
    \[
    \End(\mathcal{E}) \cong \End(\mathcal{F}) \times_{\Hom(\mathcal{F}, \mathcal{G}[1])} \End(\mathcal{G}).
    \]
    The resulting algebra is generated by the elements $t := (t, 0), \beta := (b, \beta), \gamma := (c, \gamma)$ in $\End(\mathcal{F}) \oplus \End(\mathcal{G})$, and these elements satisfy the relations of the contraction algebra.
\end{proof}

\begin{lemma}\label{lemma:Phi}
    The map $\Phi : \db(\mathcal{F} \oplus \mathcal{G}) \to \db(\Acon)$ maps $\mathcal{F}$ to the $\Acon$ module $\Acon/[\beta, \gamma]$ and $\mathcal{G}$ to the $\Acon$ module $\Acon/t$.
\end{lemma}
\begin{proof}
    Note that as a modules over $\Acon \cong \End(\mathcal{E})$, we have identifications $\End(\mathcal{F}) \cong \Acon/[\beta, \gamma]$ and $\End(\mathcal{G}) \cong \Acon/t$.

    Applying $\operatorname{Ext}^{\bullet}(-, \mathcal{F})$ to the triangle $\mathcal{E} \to \mathcal{F} \to \mathcal{G}[1]$ and using the semiorthogonality of $\mathcal{F}$ and $\mathcal{G}$, we see that the morphisms are concentrated in degree zero and there is an isomorphism of $\Acon$ modules $\Hom(\mathcal{E}, \mathcal{F}) \cong \End(\mathcal{F})$.
    Applying $\operatorname{Ext}^{\bullet}(\mathcal{E}, -)$ to the same triangle, we obtain an exact sequence
    \[
    0 \to \Hom(\mathcal{E}, \mathcal{G}) \to \End(\mathcal{E}) \to \Hom(\mathcal{E}, \mathcal{F}) \to \Ext^1(\mathcal{E}, \mathcal{G}) \to 0.
    \]
    The map $\End(\mathcal{E}) \to \Hom(\mathcal{E}, \mathcal{F})$ is the surjection $\Acon \to \Acon/[\beta, \gamma]$.

    Next, consider that $\Acon / t$ and $\Acon / [\beta, \gamma]$ are domains; in the latter case the resulting ring is a polynomial ring, while in the former case it follows from the fact that $\Acon / t$ is an Ore extension of a domain, hence is also a domain.
    It follows that the kernel for $t : \Acon \to \Acon$ is given by the ideal generated by $[\beta, \gamma]$, and the kernel for $[\beta, \gamma] : \Acon \to \Acon$ is given by the ideal generated by $t$.
    Then, the free resolution for $\Acon/[\beta, \gamma]$ is given by
    \[
    \cdots \to \Acon \xrightarrow[]{[\beta, \gamma]} \Acon \xrightarrow[]{t} \Acon \xrightarrow[]{[\beta, \gamma]} \Acon \to \Acon/[\beta, \gamma] \to 0
    \]
    and the free resolution for $\Acon/t$ is given by
    \[
    \cdots \to \Acon \xrightarrow[]{t} \Acon \xrightarrow[]{[\beta, \gamma]} \Acon \xrightarrow[]{t} \Acon \to \Acon/t \to 0
    \]
    so we see that there is a short exact sequence
    \[
    0 \to \Acon / t \to \Acon \to \Acon/[\beta, \gamma] \to 0.
    \]
    Hence, there is an isomorphism of $\Acon$ modules $\Hom(\mathcal{E}, \mathcal{G}) \cong \Acon/t$.
\end{proof}

Next, we show that this mapping induces a factorisation of the spherical functor through the map $\db(\Acon) \to \db(X_+)$.

\begin{prop} \label{prop:factor-acon}
    Let $F : \db(\mathcal{F} \oplus \mathcal{G}) \to \db(X_+) \to \db(\Lambda_0)$, and let $H : \db(\Acon) \to \db(\Lambda_0)$ be the functor induced by restriction of scalars along the morphism $\Lambda_0 \to \Acon$.
    Then the triangle below commmutes.
    \[\begin{tikzcd}
	{\db(\mathcal{F} \oplus \mathcal{G})} && {\db(\Lambda_0)} \\
	& {\db(\Acon)}
	\arrow["F", from=1-1, to=1-3]
	\arrow["\Phi"', from=1-1, to=2-2]
	\arrow["H"', from=2-2, to=1-3]
\end{tikzcd}\]
\end{prop}
\begin{proof}
    We compute the composition
    \[
    \db(\mathcal{F} \oplus \mathcal{G}) \to \db(X_+) \to \db(\Lambda_0).
    \]

    Firstly, we compute
    \begin{align*}
        \Hom_{X_+}(\sheaf \oplus V, \sheaf_{E_1}(V)) \cong \Hom_{E_1}(\sheaf \oplus V, V).
    \end{align*}
    We note that $E_1$ is the GIT quotient of the stack $S_0 \cong [V^* \oplus \Sym^2V(-1)^{\oplus 2} / GL(V)]$.
    As justified in the proof of \Cref{prop:afibfunctor}, we may compute morphisms on the stack $S_0$.
    Hence, $\Hom_{E_1}(V, V) \cong \End(\mathcal{G})$ as algebras.
    Furthermore, one sees by computing invariants that $\Hom_{E_1}(\sheaf, V) \cong \Hom_{S_0}(\sheaf, V) = 0$.
    Hence the resulting $\Lambda_0$ module structure is induced by an $\Acon$ module structure.
    One sees that the endomorphisms $t, \beta, \gamma$ of $V$ act as $0, \beta, \gamma$ respectively on $\End(\mathcal{G})$ since $t = 0$ along the locus $E_1$, so the $\Acon$ module structure is that of $\Acon/t$.
    
    Next, we compute
    \[
    \Hom_{X_+}(\sheaf \oplus V, \mathcal{F}\vert_{X_+}).
    \]
    Consider the pullback of the following diagram.
    \[\begin{tikzcd}
	& {\widetilde{S}_1} \\
	{E_2} & S
	\arrow[from=1-2, to=2-2]
	\arrow[hook, from=2-1, to=2-2]
\end{tikzcd}\]
    Here, $E_2$ is the intersection of $S$ with the semistable locus, and to take the pullback we add the conditions that $\alpha^* \neq 0$ and at least one of $\beta, \gamma$ does not preserve $\ker \alpha^*$ to the space $\widetilde{S}_1$.
    To find the resulting space, note that $\widetilde{S}_1$ is a vector bundle over $\mathbb{P}(V)$, and so the $GL(V)$ action on $\widetilde{S}_1$ can be broken to an action of a parabolic subgroup by fixing the coordinate on $\mathbb{P}(V)$, or equivalently, choosing a line in $V$.
    Since $\beta, \gamma$ preserve the associated line in $V$, upon adding the conditions we may fix $\alpha^*$ such that its kernel is a distinct line.
    In this basis, we can then write
    \begin{align*}
        & \beta = \begin{bmatrix}
            b_{00} & b_{01} \\ 0 & -b_{00}
        \end{bmatrix} & 
        \gamma = \begin{bmatrix}
            c_{00} & c_{01} \\ 0 & -c_{00}
        \end{bmatrix}
    \end{align*}
    Hence, we obtain a stack
    \[
    [\Aff^2_{b_{00}, c_{00}} \times \Aff^2_{b_{01}, c_{01}} \times \Aff^1_{p} / \CC^*]
    \]
    Note that there is a residual $\CC^*$ action coming from the subgroup $\operatorname{diag}(t, 1)$, since we still have freedom to scale the line that was chosen in $V$, and the weight matrix is given as follows.
    \[
    \begin{matrix}
                & b_{00} & c_{00} & b_{01} & c_{01} & p \\
        \lambda & 0 & 0 & -1 & -1 & 0 \\
    \end{matrix}
    \]

    The condition that $\beta, \gamma$ do not preserve $\ker \alpha^*$ becomes that $b_{01}, c_{01}$ are not both zero, so we obtain as the quotient $\mathbb{P}^1 \times \Aff^3_{p, b_{00}, c_{00}} =: \tilde{E}_2$.

    The outermost rectangle in the following diagram is Cartesian.
    \[\begin{tikzcd}
	{\tilde{E}_2} & {\widetilde{S}_1} \\
	{E_2} & S_1 \\
	{X_+} & {\thestack}
	\arrow[from=1-1, to=1-2]
	\arrow[from=1-1, to=2-1]
	\arrow[from=1-2, to=2-2]
	\arrow[from=2-1, to=2-2]
	\arrow[from=2-1, to=3-1];
	\arrow[from=2-2, to=3-2]
	\arrow[from=3-1, to=3-2]
\end{tikzcd}\]
    Applying flat base change, we find that $\mathcal{F}\vert_{X_+}$ is given by restricting $Q$ to $\tilde{E}_2$ and then pushing forward along the map $\tilde{E}_2 \to E_2 \to X_+$.
    Embedding the $\CC^*$ into $GL(V)$ via 
    \[
    \lambda \mapsto \begin{bmatrix}
        1 & 0 \\ 0 & \lambda^{-1}
    \end{bmatrix}
    \]
    we find that $Q$ restricts to $\tilde{E}_2$ as the $\sheaf_{\mathbb{P}^1}(-1)$ line bundle while $V$ restricts to $\tilde{E}_2$ as $\sheaf_{\mathbb{P}^1} \oplus \sheaf_{\mathbb{P}^1}(-1)$.

    It follows that there is an isomorphism of algebras
    \begin{align*}
        \Hom_{X_+}(\sheaf \oplus V,\mathcal{F}\vert_{X_+}) & \cong \Hom_{\tilde{E}_2}(\sheaf_{\mathbb{P}^1} \oplus (\sheaf_{\mathbb{P}^1} \oplus \sheaf_{\mathbb{P}^1}(-1)),\sheaf_{\mathbb{P}^1}(-1)) \\
        & \cong \Gamma(\tilde{E}_2, \sheaf_{\mathbb{P}^1}) \\
        & \cong \mathbb{C}[p, b_{00}, c_{00}].
    \end{align*}
    As before, the $\Lambda_0$ module structure is induced by an $\Acon$ module structure. 
    To compute the action of the endomorphisms $t, \beta, \gamma$, we note that since $t := AA^*$, it acts as $p$, and since $\beta, \gamma$ give elements of $\End_{\tilde{E}_2}(\sheaf_{\mathbb{P}^1} \oplus \sheaf_{\mathbb{P}^1}(-1))$ represented by upper triangular matrices, they act as ${b}_{00}, c_{00}$ respectively.
    Hence we identify it with the $\Acon$ module $\Acon/[\beta, \gamma]$.

    To check that the effect of the two functors agrees on morphisms, we note that it suffices to track the object $\mathcal{E}$, which is clearly sent to the $\Lambda_0$ module $\Acon$ under both functors. 
\end{proof}

Finally, we are able to show that the window-shift autoequivalence can be described as a twist around the functor from $\db(\Acon)$.

\begin{prop}\label{prop:equivalent-twists}
    The autoequivalences induced by the spherical functors $F$ and $H$ are naturally isomorphic.
\end{prop}
\begin{proof}
    The precise statement we wish to show is that there is a natural isomorphism $FF^L \simeq HH^L$, since then the dual twists induced by the spherical functors are equivalent.
    Note that $H^L = \Acon \otimes_{\Lambda_0} -$, so it defines a mapping $H^L : \db(\Lambda_0) \to \perf(\Acon)$.
    Then, $\Phi^L = \mathcal{E} \otimes_{\Acon} -$ defines a functor $\Phi^L : \perf(\Acon) \to \db(\mathcal{F} \oplus \mathcal{G})$.
    Hence, we can write
    \[
    FF^L \simeq H\Phi\Phi^LH^L
    \]
    and check that $\Phi\Phi^L \simeq \operatorname{id}$ on $\perf(\Acon)$.
    But for this it suffices to check on $\Acon \in \perf(\Acon)$, and we see that $\Phi\Phi^L(\Acon) = \Phi(\mathcal{E}) = \Acon$.
\end{proof}

\begin{remark}
In the introduction we claimed that the spherical functor we construct has some kernel, and that quotienting by this kernel allows us to recover the contraction algebra; we may now elaborate on this idea.
Consider the sheaves 
\begin{align*}
    \mathcal{F}/t & := \operatorname{coker} (t : \mathcal{F} \to \mathcal{F}) \\ 
    \mathcal{G}/[\beta, \gamma] & := \operatorname{coker}([\beta, \gamma] : \mathcal{G} \to \mathcal{G}).
\end{align*}
These two sheaves are isomorphic on the semistable locus, but differ on the unstable locus - to see this note that $\mathcal{G}/[\beta, \gamma]$ only has support contained in $\{\alpha = 0\}$, but $\mathcal{F}/t$ is also supported partly on $\{\alpha^* = 0\}$. 
It follows that the kernel of our functor $\db(\mathcal{F} \oplus \mathcal{G}) \to \db(X_+)$ contains $\operatorname{Cone}(\mathcal{F}/t \to \mathcal{G}/[\beta, \gamma])$.
Morally, we can quotient by this kernel category to obtain a smaller source category, which should then coincide with $\db(\Acon)$; as evidence, one can compute that there are no morphisms from $\mathcal{E}$ to the object given by $\operatorname{Cone}(\mathcal{F}/t \to \mathcal{G}/[\beta, \gamma])$.
\end{remark}

\appendix

\section{Locally free resolutions}

We continue to use the notation introduced in \Cref{subsubsec:strata}.

\begin{lemma}\label{lemma:resG}
    The sheaf $\mathcal{G} := \sheaf_{S_0}(V)$ has a locally free resolution:
\[
0 \to V^* \to \sheaf \oplus \Sym^2V(-1) \to V \to \sheaf_{S_0}(V) \to 0.
\]
\end{lemma}
\begin{proof}
    There is a Koszul resolution for $\sheaf_{S_0}$ given by:
\[
0 \to \sheaf(-1) \to V^* \to \sheaf \to \sheaf_{S_0} \to 0
\]
which we tensor with $V$.
\end{proof}

\begin{lemma}\label{lemma:resF}
    The sheaf $\mathcal{F} = \sigma_{*}\sheaf_{\widetilde{S}_1}(Q)$ has a locally free resolution:
    \begin{equation*} \label{eq:rewritesym}
0 \to \Sym^2V(-1) \to \sheaf^{\oplus 3} \oplus V \to V \to \sigma_{*}\sheaf_{\widetilde{S}_1}(Q) \to 0.
\end{equation*}
\end{lemma}
\begin{proof}
Firstly, we write the Koszul resolution of $\sheaf_{\widetilde{S}_1}$. 
To enforce that $\alpha$ is a subset of the line $L \subset V$, we require a vanishing of a morphism 
\[s_\alpha : \sheaf \to V \to Q.\]
Next, we require that $\beta$ and $\gamma$ preserve the line, so that we have the vanishing of two maps
\[s_{\beta}, s_{\gamma} : L \to V \to V \to Q.\]
We therefore find that $\widetilde{S}_1$ is cut out by a section of $Q \oplus (Q^2 \otimes D^{-1})^{\oplus 2}$. 
Hence, the Koszul resolution is given by:
\[
0 \to (Q^{-5} \otimes D^2) \to (Q^{-3} \otimes D)^{\oplus 2} \oplus (Q^{-4} \otimes D^2) \to (Q^{-2} \otimes D)^{\oplus 2} \oplus Q^{-1} \to \sheaf \to \sheaf_{\widetilde{S}_1} \to 0.
\]
Twisting by $Q$, we obtain:
\[
0 \to (Q^{-4} \otimes D^2) \to (Q^{-2} \otimes D)^{\oplus 2} \oplus (Q^{-3} \otimes D^2) \to (Q^{-1} \otimes D)^{\oplus 2} \oplus \sheaf \to Q \to \sheaf_{\widetilde{S}_1}(Q) \to 0.
\]
Finally, pushing down to $\thestack$ gives the exact sequence claimed, where we use the Serre duality statement $H^1(\thestack \times \mathbb{P}(V), Q^{-i}) \cong H^0(\thestack \times \mathbb{P}(V), Q^{i - 2} \otimes D)^*$.
\end{proof}

\renewcommand*{\bibfont}{\small}
\printbibliography


\end{document}